\documentclass[a4paper, 
10pt
]{article}

\usepackage[a4paper,left=3cm,right=3cm,top=2.5cm,bottom=2.5cm]{geometry}
\usepackage{amssymb}
\usepackage{amsmath}
\usepackage{graphicx}
\usepackage{subcaption}
\usepackage{dirtytalk}
\usepackage{float}
\usepackage{fancyhdr}
\usepackage{hyperref}
\usepackage{authblk}
\usepackage{tikz}
\usepackage{tikzsymbols}
\usetikzlibrary{calc,trees,positioning,arrows,chains,shapes.geometric,%
    decorations.pathreplacing,decorations.pathmorphing,shapes,%
    matrix,shapes.symbols, decorations.markings, patterns,fit}
\usepackage{color, colortbl}
\usepackage{pgfplots}
\usepackage{pgfplotstable}
\pgfplotsset{compat = 1.3}
\usepackage[ruled,vlined]{algorithm2e}

\usepackage{ulem}

\providecommand{\keywords}[1]
{
  \small    
  \textbf{\textit{Keywords---}} #1
}

\graphicspath{ {./figures/} }

\begin{document}


\title{The Neural Network shifted-Proper Orthogonal Decomposition: a Machine Learning Approach for Non-linear Reduction of Hyperbolic Equations}

\author[1,2]{Davide Papapicco}
\author[1]{Nicola Demo}
\author[1]{Michele Girfoglio}
\author[1]{Giovanni Stabile}
\author[1]{Gianluigi Rozza}

\affil[1]{Mathematics Area, mathLab, SISSA, Via Bonomea, 265, Trieste, 34136, Italy}
\affil[2]{Department of Electronics and Telecommunications, Politecnico di Torino, C.so Duca degli Abruzzi, 24, Torino, 10129, Italy}
\date{}
\maketitle

\begin{abstract}
Models with dominant advection always posed a difficult challenge for
projection-based reduced order modelling. Many methodologies that have recently
been proposed are based on the pre-processing of the full-order solutions to
accelerate the Kolmogorov $N-$width decay thereby obtaining smaller linear subspaces
with improved accuracy. These methods however must rely on the knowledge of the
characteristic speeds in phase space of the solution, limiting their range of
applicability to problems with explicit functional form for the advection field.
In this work we approach the problem of automatically detecting the correct pre-processing transformation in a statistical learning framework by implementing a deep-learning architecture. The purely data-driven method allowed us to generalise the existing approaches of linear subspace manipulation to non-linear hyperbolic problems with unknown advection fields. The proposed
algorithm has been validated against simple test cases to benchmark its performances and later successfully applied to a multiphase simulation.
\end{abstract}





\section{Introduction}\label{Intro} Reduced Order Modelling (ROM) is a
well-established set of different numerical techniques whose objective is that
of retrieving a low-rank representation of parametric differential models, s.a. Ordinary Differential Equations (ODEs) and Partial Differential Equations (PDEs), describing a relevant majority of models in physics and engineering
\cite{quarteroni_reduced_2014,stein_model_2017}. Lowering the computational cost
of numerical simulations is a fundamental aspect of both industrial and academic
research activities and therefore ROM techniques have always been, since their
introduction, an important part of the modelling process. In its most general
formulation, ROM deals with Initial Boundary Value Problems (IBVPs) in which a
PDE is parametrised by a finite set of $P\in\mathbb{N}$ parameters. One such
example is the following scalar, linear, first-order in space IBVP with
non-parametric, steady Dirichlet's Boundary Conditions (BCs) and Initial
Condition (IC)
\begin{equation}\label{eq1}
\begin{cases}
    \partial_t u(\mathbf{x},t,\boldsymbol{\mu}) + \mathcal{L}(u,\boldsymbol{\mu}) = 0\,,\quad\mathbf{x}\in\Omega\subset\mathbb{R}^d\,,\:t\in[0,T>0]\,,\:\boldsymbol{\mu}\in\mathcal{P}\,, \\
    u(\mathbf{x},t,\boldsymbol{\mu}) = g(\mathbf{x})\,,\quad\forall\mathbf{x}\in\partial\Omega\,,\\
    u(\mathbf{x},t=0,\boldsymbol{\mu}) = u_0(\mathbf{x})\,,\quad\forall\mathbf{x}\in\Omega\,,
\end{cases}
\end{equation}
where $d=1,2,3$ represents the number of spatial dimensions of the model, $T>0$
is a final instant that identifies the interval in which the time evolution of
the model is evaluated and $\mathcal{P}\subset \mathbb{R}^P$ is the parameter domain of the PDE s.t. $\dim(\mathcal{P})=P$. Following a domain discretation of $\Omega$ using e.g. Finite Volumes, Finite Elements or Finite Differences, one obtains a numerical
approximation $\mathbf{u}_h(t,\boldsymbol{\mu})\in\mathcal{V}_h$ for the
(discrete) time evolution of $u(\mathbf{x},t,\boldsymbol{\mu})$ as evaluated at
timesteps $t_k\,,\;k=0,1,\dots,$ and where $N_h := \dim(\mathcal{V}_h) $
represents the Degrees of Freedom (DOFs) of the problem. By defining the discrete
solution manifold $\mathcal{M}_h := \big\{\mathbf{u}_h(t_k,\boldsymbol{\mu})\in\mathcal{V}_h\,,\:k\in\mathbb{N}_0\,,\:\;\text{s.t.}\:\,\boldsymbol{\mu}\in\mathcal{P}\big\}$ one immediately recognises that a high-dimensional approximation space $\mathcal{V}_h$ leads to an expensive computational cost whenever one wants to evaluate the time evolution of $\mathbf{u}_h(t,\boldsymbol{\mu})$ for multiple
instances of $\boldsymbol{\mu}$ in the parameter space $\mathcal{P}$.
Formally stated, ROM aims at reducing the computational cost of a (numerical)
parametric simulation by constructing the best low-rank approximation
$\mathcal{R}$ of manifold $\mathcal{M}$ s.t. it encodes the parameters variation
of the Full Order Model (FOM) but with a reduced number of DOFs. One such
example is that of evaluating the pressure field solution in response to several
geometrical configurations of an airfoil parametrised by e.g. different angle of
attack, chord line and mean thickness. Many different classes of ROM
techniques have been developed for the low-rank reduction of several phenomena
\cite{quarteroni_reduced_2016,hesthaven_certified_2016,stabile_pod-galerkin_2017}
and in particular projection based ROM have found widespread application in
reducing different parametric differential problems modelling several phenomena in fluid dynamics. This class is based on the identification of a reduced set of basis functions, or modes,
$\big\{\boldsymbol{\phi}_j\big\}_{j=1,\dots,R}$ s.t. their superposition spans
the best possible low-rank approximation of the solution manifold i.e.
$\text{span}(\boldsymbol{\phi}_1,\dots,\boldsymbol{\phi}_R) = \mathcal{R}
\approx\mathcal{M}_h$. The Proper Orthogonal Decomposition (POD), on which this
work is based upon, is one of those techniques in which the
basis extraction is performed using Singular Value Decompostion (SVD) of a
so-called snapshot matrix $\boldsymbol{\mathcal{X}}$ in which the FOM solutions
$\mathbf{u}_h(t,\boldsymbol{\mu})$ are stored as column vectors. When paired
with a projection technique
$\Pi_{\mathcal{R}}\mathbf{u}_h(t,\boldsymbol{\mu})\in\mathcal{R}$ s.a. the
Galerkin projection the resulting POD-Galerkin method allows for fast
computation of several simulations with different instances of
$\boldsymbol{\mu}$ in the domain $\mathcal{P}$. The efficiency of any
ROM technique, including POD-Galerkin methods, stems directly from the number of
basis functions needed to encode the parametric essential dynamics of a model at
reduced order. The reduction accuracy of a model has been traditionally
quantified by the Kolmogorov $N-$width of the problem which is defined, in
$L_2-$norm, as
\begin{equation}\label{eq2}
    \text{dist}_R(\mathcal{M}_h):=\inf_{\mathcal{R}\in\mathcal{R}_R}\sqrt{\frac{\sum_{\mathbf{u}_h\in\mathcal{M}_h}||\mathbf{u}_h-\Pi_{\mathcal{R}}\mathbf{u}_h||^2}{\sum_{\mathbf{u}_h\in\mathcal{M}_h}||\mathbf{u}_h||^2}}\,.
\end{equation}
Problems that are characterised by a fast Kolmogorov $N-$width decay (here we used $N\equiv R$ to be in line with the notation commonly found in the literature) are easily
restricted by a low-rank linear subspace approximation; on the other hand a slow
decay indicates that a linear approximation assumption might be inappropriate
and inaccurate for that problem. For this reason, despite the popularity that
POD-Galerkin methods gained in various areas of computational science and
engineering
\cite{girfoglio_pod-galerkin_2020,ballarin_supremizer_2015,ballarin_pod-galerkin_2016,hijazi_data-driven_2020},
they have been severely limited by the fact that the low-rank approximation of
$\mathcal{M}_h$ they build is sought within a sequence
$\big\{\mathcal{R}_R\big\}_{R=1,\dots,\infty}$ of $R-$dimensional linear
subspaces. The assumption that a FOM can be accurately approximated by
restricting its dynamics on a low-rank linear subspace proved in fact to be
suitable for elliptic and parabolic problems (where the Kolmogorov $N-$width decays
exponentially with $R$ \cite{cohen_kolmogorov_2015}) while it is not trivially
extendable to problems with dominant advection s.a. those modelled by hyperbolic
PDEs (for which the Kolmogorov $N-$width decay is much less steep
\cite{nonino_overcoming_2019}). It results that the reduction of
advection-dominated PDEs with POD-Galerkin methods requires a large number $R$
of basis functions to be accurately depicted as they struggle to restrict the
FOM to a linear subspace approximation of $\mathcal{M}_h$, essentially nullifying
the reduction in computational cost. Since models with dominant advection are
particularly recurrent in fluid dynamics, being associated to the conservation
laws describing a large multitude of phenomena (most notably Euler's equations
and Riemann's problems, shallow water equations, multiphase models), a
growing number of endeavours in recent years have proposed alternative
modifications for the improvement in performances for the approximation of such
manifolds
\cite{iollo_advection_2014,pacciarini_reduced_2015,torlo_stabilized_2018}. The
work that has been done is very heterogeneous in terms of the methodologies
adopted, with many efforts that brought in the latest data-driven trends in
machine learning and features extraction. Nevertheless the state-of-the-art in
improving the accuracy of ROM of hyperbolic equations can be distinct into two
major approaches: those based on ad-hoc transformations of the linear subspace
as a support for better mode extraction of projection-based ROM
\cite{rim_model_2018,rim_transport_2018,chetverushkin_model_2019,nair_transported_2019,taddei_registration_2020}
and those based on a construction of non-linear manifolds
\cite{kashima_non-linear_2016,hartman_deep-learning_2017,crisovan_model_2018,hoang_projection_2021}.
The former class has been the
first to be approached while the latter strategy is only most recently being
explored to overcome the limitations of the earlier approaches. Among the
various reasons there is the fact that the underlying methodology of
pre-processing the linear subspace has been that of embedding the dominant
advection of the solution at offline stage using a deterministic transformation.
This requires some sort of problem-specific prior knowledge of the physical or
mathematical properties of the FOM e.g. the characteristic speeds in phase space
of a conservation law \cite{reiss_shifted_2018,sarna_hyper-reduction_2021}. This
obviously limits the range of applicability of pre-processing based model
reduction to deal with hyperbolic PDEs that feature simple advection fields. To
overcome such shortcomings the second class of methods based on deriving
non-linear manifold reduction directly has been recently favoured with important
results being obtained with the adoption of convolutional autoencoders
\cite{lee_model_2019} and Arbitrary Lagrangian-Eulerian (ALEs) frameworks
\cite{torlo_model_2020,mojgani_arbitrary_2017}.
  In this work we introduce a novel statistical learning framework for
the generalisation of the earliest pre-processing techniques to those hyperbolic
problems that feature non-linear and unknown advection fields. The resulted
Neural Network shifted-Proper Orthogonal Decomposition (NNsPOD) algorithm is a
fully non-intrusive, purely data-driven machine learning method that seeks for
an optimal mapping of the various snapshots in the low-rank linear subspace to a
reference configuration via an automatic detection that does not depend on the
physical properties of the model. By construction our algorithm does not belong
exclusively to either of the aforementioned classes of manifold reduction. While
NNsPOD starts with the goal of building a pre-processing transformation that
accelerates the Kolmogorov $R-$width decay for advection-dominated problems, it
also achieves a virtually unlimited range of applicability to non-linear
manifold reduction since it requires no prior knowledge on the properties nor
complexity of the FOM.
During the development of our algorithm we became aware of a project being proposed in \cite{peng_learning-based_2021} which achieves the similar result of non-linear manifold reconstruction of transport-dominated problems via deep-learning models. Despite being both projection-based and physics-independent, our methodology substantially differs from the one proposed in \cite{peng_learning-based_2021} in that NNsPOD uses a neural architecture to learn an optimal pre-processing transformation for the snapshots. On the other hand, as highlighted in Section 3 of \cite{peng_learning-based_2021}, their algorithm focuses on using two neural architectures to learn adaptively the parametric reduced basis of the FOM at offline stage and then combine their outputs through an inner product to obtain the ROM. We also emphasize that all the benchmarks that follow in the present work are $2-$dimensional problems as opposed to the $1-$dimensional problems presented in \cite{peng_learning-based_2021}.
The present work is structured as follows: in Section~\ref{Pre-processing} we contextualise the state-of-the-art of linear subspaces transformations while also showing the limitations of traditional POD-Galerkin reduction for a canonical problem with dominant advection; in Section~\ref{NNsPOD} we present the methodologies adopted by NNsPOD as a generalisation of such techniques to unknown and non-linear transport fields; in Section~\ref{Multiphase} we demonstrate the performances of NNsPOD as a result of its application to a reduction of a multiphase model; in Section~\ref{Conclusions} we conclude by highlighting the outstanding challenges of our approach while also providing some insights on practical limitations we encountered during the development and test stage of our algorithm.

\section{Pre-processing transformation of POD linear
subspaces}\label{Pre-processing} In order to introduce the new methodology
proposed in this work we shall contextualise the common approach shared by all
the previous works that have been done in linear subspace pre-processing of
advection-dominated models. To that end we start by formally introducing the
standard procedure of a traditional POD-Galerkin technique and highlight its
under-performance with a simple hyperbolic equation.
\subsection{Full-order scalar advection equation}\label{sub2.1}
Let us consider a simple test case of the IBVP in $(\ref{eq1})$ for which
$\mathcal{L} = \nabla\cdot\mathbf{F}$ i.e.
\begin{equation}\label{eq3}
    \begin{cases}
        \partial_t u(\mathbf{x},t) + \nabla\cdot\mathbf{F} = 0 \,,\quad\mathbf{x}\in[0,1]\times[0,1]\equiv\Omega\,,\:t\in[0,T=1]\,,\\
        u(\mathbf{x},t) = 0\,,\quad\forall\mathbf{x}\in\partial\Omega\,,\\
        u(\mathbf{x},0) = \text{exp}(-\frac{1}{2}\,\mathbf{x}^T\mathbf{x})\,,\quad\forall\mathbf{x}\in\Gamma\subset\Omega\,,
\end{cases}
\end{equation}
which is a linear first-order advection equation in a scalar unknown
$u(\mathbf{x},t)$ that is null everywhere over a unitary square domain except for a bounded region $\Gamma$ in which varies as a multivariate Gaussian pulse. 
\begin{figure}[H]
    \centering
    \includegraphics[keepaspectratio,width=\textwidth]{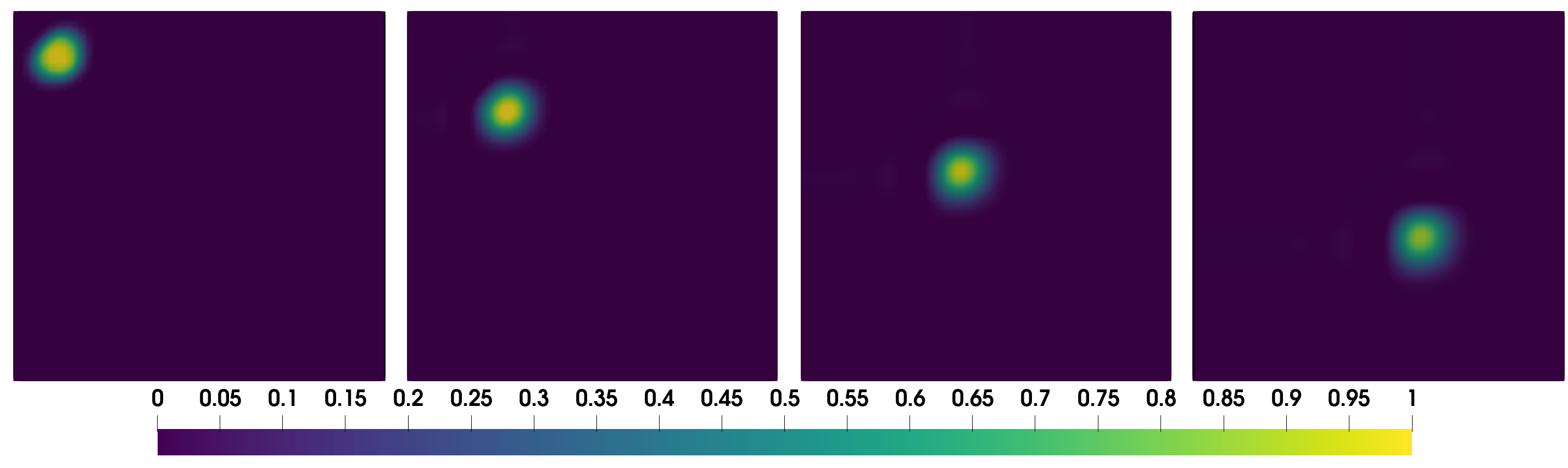}
    \caption{Different snapshots of the FOM simulation in $(\ref{eq3})$ sorted left to right from         the IC to increasing timesteps.}
    \label{fig1}
\end{figure}
For the sake of simplicity we set $(\ref{eq3})$ to have only one parameter which is the time
variable itself; this entails that $P=dim(\mathcal{M}_h)=1$ which in turn means
that we expect to retrieve a sufficiently accurate linear subspace approximation
of the FOM with only $R=1$ reduced basis vector. We know however that this is
not the case since $(\ref{eq3})$ is an hyperbolic PDE.
To test the inability of traditional projection-based ROM at
restricting the FOM to a linear subspace let us set a computational grid of
$50\times50$ points in $\Omega$ while we discretise the time interval into $100$
uniform steps to comply with the Courant –Friedrichs–Lewy (CFL) condition. We
also set the vector field $\mathbf{F} = \mathbf{f}\,u$ to be a linear flux
function that is uniform in space and constant in time with the advection field
being $\mathbf{f} = (1,-1)$. We obtain a FOM numerical solution
$\mathbf{u}_h(t)\in\mathcal{V}_h\,,\;N_h = 2500$ using a Finite Volume (FV)
spatial discretisation scheme, while the time evolution over the discrete
timesteps is obtained via an unconditionally stable, first-order accurate,
implicit Euler method. For the sake of simplicity we will not discuss the
theoretical background of both numerical methods as it is out of the scope of
the present work; the reader may refer to
\cite{moukalled_finite_2016,atkinson_num_anal,quarteroni_numerical_2009} for
a comprehensive reading on standard algorithms in introductory and advanced
numerical analysis. In order to present a particularly challenging advection
model to our POD-Galerkin reduction method, we reduce as much as possible the
numerical diffusion given by low-order cell face interpolation; as such a
third-order accurate QUICK scheme \cite{leonard_quick_1979} is adopted with the
numerical stability being assured by the combined implicit Euler time advancing
method and CFL-complying timesteps. Different FOM snapshots of this simulation
are depicted in Figure~\ref{fig1}.
\subsection{The POD method for reduced basis extraction}\label{sub2.2}
The standard procedure of POD reduction is hereby outlined. The $N_s=100$ snapshots, sampled in the $1-$dimensional FOM manifold during the simulation, are sorted by timestep in increasing order and stored as column vectors of a $2500\times100$ snapshot matrix. We can find the best low-rank approximation of $\boldsymbol{\mathcal{X}}$ via SVD as per the Eckart–Young theorem \cite{stewart_svd_1993}
\begin{equation}\label{eq4}
    \boldsymbol{\mathcal{X}} = \big(\mathbf{u}_h(t_1=0),\dots,\mathbf{u}_h(t_{100}=T)\big) \approx \mathbf{L}\boldsymbol{\Sigma}\mathbf{R}\,,\;\mathbf{L}\in\mathbb{R}^{N_h\times R}\,,\:R=1,\dots,N_s.
\end{equation}

\begin{figure}[H]
    \centering
    \includegraphics[keepaspectratio,width=\textwidth]{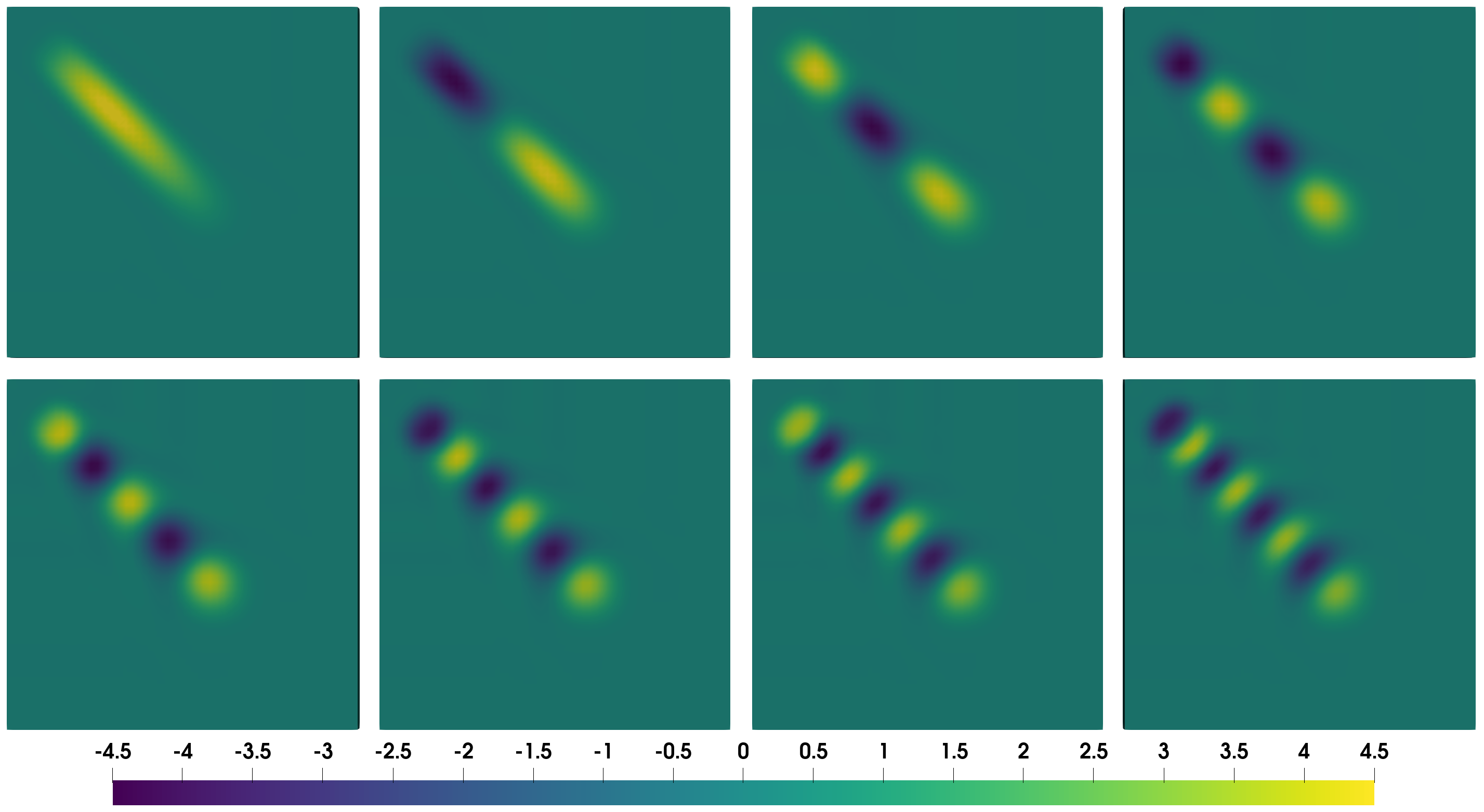}
    \caption{First 8 modes extracted via POD reduction from the FOM in $(\ref{eq3})$ sorted in ascending order from top to bottom and from left to right.}
    \label{fig2}
\end{figure}

The left singular vectors in $\mathbf{L}$ are those that are associated with any subspace approximation $\mathcal{R}_R$ of $\mathcal{M}_h$; in Figure~\ref{fig2} we depict the first $8$ of those. As specified in~\ref{sub2.1}, a linear combination of $R=P=1$ of those vectors should, theoretically, retrieve a reduced order model with high accuracy. However if we set such threshold to be e.g. $10^{-3}$ we immediately realise that such degree of accuracy is obtained with a reduced basis of $14$ left singular vectors. This model clearly cannot be constrained easily in a linear subspace as produced by a traditional POD method. 

\subsection{Linear subspace shift-based transformation}\label{sub2.3}
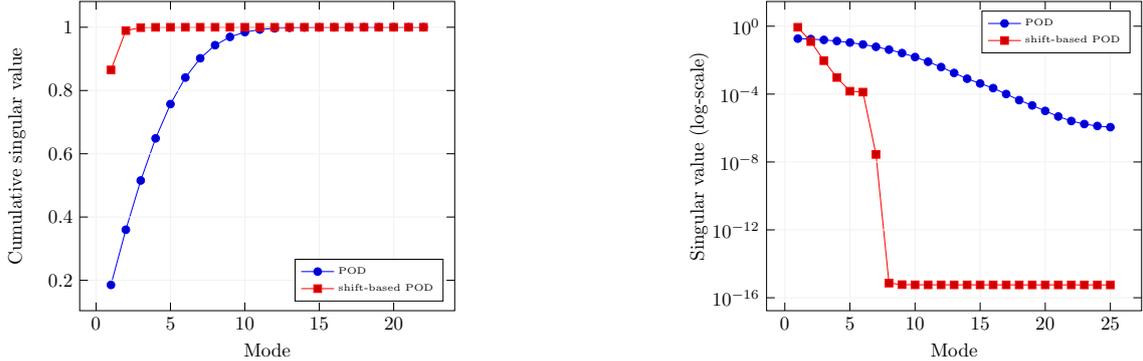
\begin{figure}[H]
\begin{tikzpicture}[scale=0.72]

\definecolor{color1}{rgb}{0.75,0,0.75}
\definecolor{color0}{rgb}{0,0.75,0.75}

\begin{axis}[
axis on top,
legend cell align={left},
legend style={font=\tiny},
tick pos=both,
xtick style={color=black},
ytick style={color=black},
xlabel={Mode},
grid=both,
legend pos=south east,
grid style={line width=.1pt, draw=gray!10},
ylabel shift = 5 pt,
ylabel={Cumulative singular value}
]
\addplot table [x expr=\coordindex+1, y index=0] {data/Square_CumSV_POD};
\addlegendentry{POD}
\addplot table [x expr=\coordindex+1, y index=0] {data/Square_CumSV_sPOD};
\addlegendentry{shift-based POD}
\end{axis}
\end{tikzpicture}\hfill
\begin{tikzpicture}[scale=0.72]

\definecolor{color1}{rgb}{0.75,0,0.75}
\definecolor{color0}{rgb}{0,0.75,0.75}

\begin{axis}[
axis on top,
legend cell align={left},
legend style={font=\tiny},
tick pos=both,
xtick style={color=black},
ytick style={color=black},
xlabel={Mode},
ymode=log,
grid=both,
legend pos=north east,
grid style={line width=.1pt, draw=gray!10},
ylabel shift = -5pt,
ylabel={Singular value (log-scale)}
]
\addplot table [x expr=\coordindex+1, y index=0] {data/Square_SV_POD};
\addlegendentry{POD}
\addplot table [x expr=\coordindex+1, y index=0] {data/Square_SV_sPOD};
\addlegendentry{shift-based POD}
\end{axis}

\end{tikzpicture}
\caption{Comparative accuracy of traditional POD and shift-based pre-processing POD for the same advection model in $(\ref{eq3})$.}
\label{fig3}
\end{figure}

In order to improve the approximation of $\mathcal{M}_h$ via linear subspaces,
different methods have been proposed in the most recent years; the majority of the
them started to focus in deriving an ad-hoc transformation of the linear
subspace in a frame of reference that facilitates the decay of the Kolmogorov
width. This shift-based pre-processing is in fact consistent with the hyperbolic
character of models with dominant advection which are solved at full-order, both
analytically and numerically, by exploiting the characteristics curves in phase
space of the solution. Therefore, as long as the characteristic speeds of an
hyperbolic PDE are known, one can pre-process the snapshots of the FOM in
$\boldsymbol{\mathcal{X}}$ by e.g. mapping them to the IC following backward
those curves in phase space. This approach therefore aims at finding a better
frame of reference for the SVD of the POD by transporting the FOM snapshots by
means of: backward-shift transformation
\cite{reiss_shifted_2018,sarna_hyper-reduction_2021} and successive
interpolation; displacement interpolation technique \cite{rim_model_2018} or
shock-fitting methods \cite{chetverushkin_model_2019}; locally adapted bases
\cite{peherstorfer_model_2020}.
  In the following we will focus on the first of those linear subspace pre-processing methods. As explained in \cite{reiss_shifted_2018} one builds a discrete shift operator $T_{\mathbf{b}}$ that acts on a space-time dependent function $u(\mathbf{x},t)$ transporting its frame of reference by an amount that is proportional to the transport field $\mathbf{b}$. The $k-$th FOM snapshot in $(\ref{eq4})$ is an $N_h-$dimensional vector; each of its component store the field value of the numerical solution associated to the $x$ and $y$ coordinates of a centroid $\mathbf{x}$ on the computational grid following the FV discretisation. Applying the shift operator to such snapshot amounts to derive a new spatial distribution of those centroids for the same field values i.e.
\begin{equation}\label{eq5}
    T_{\mathbf{b}}\,\mathbf{u}_h(t_k) = \tilde{\mathbf{u}}_h = \big(u_h(\mathbf{x}_j - \mathbf{b}t_k,0)\big)_{j=1,\dots,N_h}\,,\quad\forall k=1,\dots,N_s\,.
\end{equation}
Furthermore an interpolation of $\tilde{\mathbf{u}}_h$ is necessary in order to
reconstructs the FOM field solution starting from the shifted points
distribution on the manifold; the transformed snapshot matrix is thus decomposed
in the new frame of reference by SVD. As depicted in Figure~\ref{fig3} in fact,
a shift transformation on the snapshots allows for a better approximation with
just $4$ modes needed to span a linear subspace that is within the same
$10^{-3}$ accuracy threshold set in~\ref{sub2.2} for $\mathcal{M}_h$.

\subsection{Limitations in deducing the shift operator}\label{sub2.4}
The stated claim in \cite{reiss_shifted_2018} is to build a data-driven,
shift-based, pre-processing transformation that identifies the dominant
structures of a hyperbolic model regardless of the characteristic velocities,
with the goal of extending the methodology to non-linear advection fields. The
authors also add that in order to identify those \textsl{ansatz} modes, whenever
the phase-space velocities are unknown explicitly, two different approaches
might be considered. One approach is based on the tracking of the peaks in the
field values of the numerical solution across different snapshots; this method
is although unfeasible for hyperbolic models that do not preserve the time or
parameter dependent shape of the IC across their time evolution (e.g. multiphase
simulations as we address in Section~\ref{Multiphase}). The other suggested
methodology, that is actually discussed in \cite{reiss_shifted_2018}, is to
perform multiple SVDs with sampled values of the shift velocities; in a
data-driven fashion one should thus be able to detect the correct backward shift
transformation by examining the singular values spectrum and isolate those for
which the decay is maximised. This approach is indeed feasible to deal with
hyperbolic PDEs that feature unknown transport velocities. However we identify
two major limitations:
\begin{itemize}
    \item The examination of the singular value spectrum as a function of trial shift operators entails an unnecessary increase in computational complexity during the offline stage that scales polynomial with $d=1,2,3$. This is due to the fact that multiple snapshot matrices, and subsequent SVD, are required in order to detect the appropriate spatial transformation that maximise the Kolmogorov width decay since, in order to achieve a truly data-driven degree of flexibility, no prior information on the physical behaviour shall be used to sample the velocities for the trial shift operators.
    \item  The various guesses for candidate shift operators are all, as stated by the authors in \cite{reiss_shifted_2018}, limited to uniform and constants advection fields, although the shift itself can easily be computed pointwise in the grid and at different timesteps. This  introduces a further degree of approximation and loss in accuracy for the case of non-linear hyperbolic PDEs; in fact guessing linear, constant and uniform shift velocities essentially amounts in linearising those characteristic curves in phase-space that may present non-trivially integrable irregularities.
\end{itemize}
Our aim is that of generalising the shift-based pre-processing linear subspace approximation by overcoming these aforementioned limitations in performance, computational scalability and models of advection. In our approach we in fact retrieve a non-linear transformation that achieves the same results in \cite{reiss_shifted_2018} but with an automatic detection that is both consistent with the method of characteristics and also robust to be applied to hyperbolic PDEs with unknown and non-linear advection.

\section{Automatic detection of bijective mapping in non-linear manifolds}\label{NNsPOD}
The traditional fields of application for statistical learning models, since their renewed development, have been data-analysis and data-mining, signal processing and computer vision. In the last years however they found widespread implementation in different subsets of scientific computing as well \cite{frank_machine-learning_2020}. Industrial-driven requirements in improving the computational cost of numerical simulations in engineering, natural and life sciences resulted in an ever increasing integration of machine-learning algorithms within the more "traditional" numerical methods. The introduction of Artificial Neural Networks (ANNs), and specifically Deep Neural Networks (DNNs), further encouraged this coupling by providing purely data-driven non-linear mappings. ROM has been no exception to this merge of disciplines and the similarities between POD and Principal Component Analysis (PCA) is one example that justifies the natural blending of one set of techniques into the other. Given the extensive presence of complex non-linear PDEs, s.a. those modelling turbulent regimes \cite{kutz_deep_2017}, Computational Fluid Dynamics (CFD) represented the quintessential field of experimentation for ANNs integration in scientific computing with the first endeavours \cite{sarghini_neural_2003} being introduced as early as 2003.
The DNN-based algorithm we hereby propose has been developed along those directions outlined above, that is by integrating, in an efficient but mathematically consistent process, a statistical learning paradigm into a traditional ROM technique s.a. POD. The main limitation in detecting the appropriate backward transformation for shift-based POD reduction of models like $(\ref{eq3})$ is that prior knowledge of the characteristic velocities is necessary in order to quantify the pre-processing map. A data-driven approach based on singular values spectrum analysis does not automate the process as, in order to be efficient, a proper sampling of the velocity space has to be performed. Furthermore if there is a parametric dependence of such transport field, i.e. $\mathbf{b} = \mathbf{b}(\boldsymbol{\mu},t)$, then the Galerkin projection will not generalise well during the online phase since the pre-processed manifold does not contain the information for the evolution of the numerical solution along different characteristic curves in phase-space. The approach taken in previous works s.a. \cite{reiss_shifted_2018,rim_transport_2018} can thus be modelled as follows
\begin{equation*}
    T_{\mathbf{b}}:\mathcal{M}_h\mapsto\tilde{\mathcal{M}}_h\,,\:\tilde{\mathcal{M}_h}:=\big\{\mathbf{u}_h\big(\mathbf{x}_j-\mathbf{b}(t_k,\boldsymbol{\mu}_p\big)\,t_k)\,,\:\forall k=1,\dots,N_s,\,p=1,\dots,N_p\big\}\,,
\end{equation*}
i.e. finding a backward map $T_{\mathbf{b}}$ in which the $N_p$ parametric instances of the transport field at offline stage is embedded. This method changes the frame of reference of the solution manifold making difficult its generalisation at online stage.
  A deep-learning framework on the other hand can indeed derive a bijective transformation between the shifted manifold (of which the linear subspace approximation is build) and the original one by ignoring such parametric dependence for the advection field $\mathbf{b}$. This motivates the needs for the adoption of a non-intrusive technique for the computation of the backward map based on automatic shift-detection i.e. a machine learning algorithm that does not quantify the shift operator neither based on the pointwise value of the transport field nor a linear approximation of it. In our approach we thus seek to derive instead a bijective mapping 
\begin{equation*}
    \mathcal{C}^1\ni\mathcal{T}:\mathcal{M}_h\mapsto\tilde{\mathcal{M}_h}\,,
\end{equation*}
which does not linearise the characteristic velocities in phase-space. In the development of NNsPOD we focused in assuming that no prior information regarding the functional dependence of $\mathbf{b}$ is known as we seek to derive a bijection that traces back any snapshot in $\boldsymbol{\mathcal{X}}$ to an arbitrary \textsl{reference configuration} $\mathbf{u}_{\text{ref}}\in\mathcal{M}_h$. It might seem trivial to choose the IC as reference configuration; however, as we will explain further on in the section, given the general framework in which NNsPOD was conceived, we will 
consider also other snapshots as reference configuration for improving the bijective map. 
Being NNsPOD a data-driven, pre-processed ROM it should also provide a consistently fast online phase for parametric formulation of time-dependent hyperbolic PDEs. We postpone the results of such models, as well as the integration of NNsPOD with a Galerkin projection, to our future works; 
in the present work we only highlight the generalisation of the effectiveness of shift-based linear approximations in POD reduction methods for particularly difficult advection-dominated models with unknown and non-linear transport fields (as presented in Section~\ref{Multiphase}).
\subsection{Statistical learning formulation and reference configuration}\label{sub3.1}
Our framework is naturally implemented as a statistical learning technique in which the objective is to detect automatically an optimal transformation of the solution manifold $\mathcal{M}_h$ that increases the Kolmogorov width decay associated to its model. If a DNN is deployed to build such bijective map while it merely relies on the snapshots collected at FOM, then the desired requirements are: 
\begin{itemize}
    \item It must preserve the hyperbolic nature of the PDE consistently with the dominant advection model: NNsPOD has to be flexible enough to output transformations that have to be, in principle, rigid transports of any collected snapshot to a reference configuration of choice without adding numerical diffusion while, at the same time, being able to also change the shape of those snapshots if it is required (e.g. multiphase simulation we will discuss in Section~\ref{Multiphase}).
    \item The data-flow for the two phases of shift-detection and field-reconstruction (interpolation), as outlined in \cite{reiss_shifted_2018} has to be continuous: in order for the backward map to be generalised for the Galerkin projection the non-linear transformation provided by NNsPOD has to be invertible for any input in the training set and as such a continuity constraint in the architecture of the DNN is imposed.
\end{itemize}
In order to encode the automatic shift-detection in the backward map itself we must convert the variational form of singular values spectrum analysis and velocity sampling
into a statistical learning theoretical framework. First and foremost we interpret the FOM snapshots as data-points of the training set for the network
\begin{equation}\label{eq6}
    \mathbf{M}:=\big(\mathbf{u}_h(\boldsymbol{\mu}_j)\big)_{j=1,\dots,P}\,,\quad\mathbf{X}:=\boldsymbol{\mathcal{X}}^T\in\mathbb{R}^{N_s\times N_h}\,.
\end{equation}
We observe that the training set is itself a subspace of the full-order solution manifold; in statistical learning theory we will thus have a $N_h-$dimensional features space of cardinality $N_s$. Secondly we assign a semi-supervised learning paradigm to our model by choosing arbitrarily, among the samples in $\mathbf{M}$, a reference configuration $\mathbf{u}_{\text{ref}}$ for the field values that acts as unique label for the training of the network on the rest of the collected data-points. It is important to clarify this aspect as it is one of the features of NNsPOD that can be exploited in order to build more refined and efficient machine learning algorithms in ROM of hyperbolic equations. 
For the sake of simplicity let us consider a $1-$dimensional formulation of the IBVP in $(\ref{eq3})$ in which the advection-dominated model features a regular, constant and uniform transport field $b=1$. Given a stable full-order numerical simulation (i.e. with very low numerical diffusivity) it is intuitive that the reference configuration and any given
data-point in $\mathbf{M}$ only differ for the centroids $x$ at which a certain value of field solution is associated. For instance, the peak of the gaussian pulse will be located, e.g. at
$x = 0$ for a reference configuration $\mathbf{u}_h(x, t_0)$ and we can easily predict that, after $\Delta t=1$, it will be located at $x = 1$ i.e. $\mathbf{u}_h(x, t_1 = t_0 + 1)$. It is trivial to see that, in this instance, choosing the IC as reference configuration is no different then choosing any other snapshots in $\mathbf{M}$ when building a backward map either by direct calculation or automatic detection. This arbitrariness is itself an advantageous property of the hyperbolicity of the equation that is easily exploited whenever there are simple transport fields as stated above. However, if one is dealing with the more concrete case of irregular transport fields, then the IC itself does not have to be necessarily the one that fully reconstructs the shifted manifold since the stationary frame of reference can be associated to any of the snapshots collected at FOM. To take full advantage of faster Kolmogorov width decay, NNsPOD does not restrict itself in sampling backward transformations to the IC but to any reference configuration that optimises the reduction; the non-linear projection operator associated to the automatically detected map will thus be able to generalise from any frame of reference thanks to the bijective constraint discussed above.
  For this reason in the following we will refer to the unique label, chosen for the training of the neural network part of NNsPOD, in more general terms as the reference configuration $\mathbf{u}_{\text{ref}}\in\mathbf{M}$ which, we reiterate, does not necessarily have to coincide with the IC of the IBVP. The setting is completed by choosing one particular metric for computing the loss function between the output associated to any data-point in $\mathbf{M}$ and the reference configuration. While there is substantial space for testing and experimenting with this choice,
in the development of the present work we restrict the development to the $L_2-$norm which is consistent with the Kolmogorov width decay of the finite-dimensional functional space $\mathcal{V}_h$ for the FOM snapshots
\begin{equation}\label{eq7}
    J(\mathbf{x}_{\text{ref}},N_s):=\frac{1}{N_s}\sum_{j=1}^{N_s-1}|| \tilde{\mathbf{u}}_h(\boldsymbol{\mu}_j) - \mathbf{u}_{\text{ref}}||_2\,,\quad \tilde{\mathbf{u}}_h := \mathbf{u}_h(\mathbf{x},\boldsymbol{\mu})\circ\mathcal{T}
\end{equation}
Future works might be able to derive a rigorous mathematical structure for the selection of appropriate metrics for the training of shift-detecting networks in projection-based ROM pre-processing.
\subsection{Architectures for continuous data-flow: the shift-detection and field-  reconstruction split}\label{3.2}
Aside from the shift detection itself, one important procedure for a pre-processing transformation of the FOM manifold is the ability of reconstructing the field values of the solution following the backward map to the reference configuration. This can be achieved via e.g. interpolation \cite{reiss_shifted_2018} of the field values from the shifted points to the nearest centroids on the grid. In developing an automatic shift-detection algorithm we therefore must include such field-reconstruction part within the algorithm itself. Being purely data-driven, it is also desirable that the interpolation is also based on the FOM snapshots and does not rely on the prior physical knowledge of the problem (e.g. upwind interpolation methods). One of the major advantages of NNsPOD is that the choice of the reference configuration for the training of the network is arbitrary. We exploit such potential by splitting the workload between two separate networks:
\begin{itemize}
    \item \textsl{ShiftNet} is the neural network that has the duty of quantifying the optimal shift for pre-processing transformation of the full-order manifold and that maximises the Kolmogorov width decay.
    \item \textsl{InterpNet} is the neural network that must \say{learn} the reference configuration in the best possible way w.r.t. its grid point distribution s.t. it will be able to reconstruct its \say{shape} for every shifted centroid distribution.
\end{itemize}
\begin{algorithm}[H]
\SetAlgoLined
\KwResult{Optimal shift-detected transformation $\mathcal{T}$}
 Construction of snapshot matrix $\boldsymbol{\mathcal{X}}$\;
 Computing the feature matrix $\mathbf{X}=\boldsymbol{\mathcal{X}}^T$\;
 Setting the accuracy thresholds $\varepsilon_{\text{SVD}}$, $\varepsilon_{\text{shift}}$, $\varepsilon_{\text{interp}}$\;
 \While{$\varepsilon>\varepsilon_{\text{SVD}}$}{
    Selection of reference configuration $\mathbf{u}_{\text{ref}}$\;
    \While{$\varepsilon^{'}>\varepsilon_{\text{interp}}$}{
        \textbf{InterpNet}.forward\;
        compute $\varepsilon^{'}$\;
        \textbf{InterpNet}.backward\;
    }
    $\mathcal{T_{\text{interp}}} = $ \textbf{InterpNet}.forward\;
    \While{$\varepsilon^{''}>\varepsilon_{\text{shift}}$}{
        $\tilde{\mathbf{x}} =$ \textbf{ShiftNet}.forward\;
        $\tilde{\mathbf{u}}_h = \mathcal{T_{\text{interp}}}\circ \tilde{\mathbf{x}}$\;
        compute $\varepsilon^{''}$\;
        \textbf{ShiftNet}.backward\;
    }
    $\mathcal{T} = $ \textbf{ShiftNet}.forward\;
  \For{$\mathbf{u}_h\in\mathbf{X}$}{
    $\tilde{\mathbf{u}}_h = \mathbf{u}_h\circ\mathcal{T}$\;
    }
    Perform the SVD on $\tilde{\mathbf{X}}$\;
    Compute $\varepsilon_{\text{SVD}}$\;
 }
 \caption{The NNsPOD algorithm}\label{alg1}
\end{algorithm}

\begin{figure}[H]
    \centering
    \includegraphics[keepaspectratio,width=\textwidth]{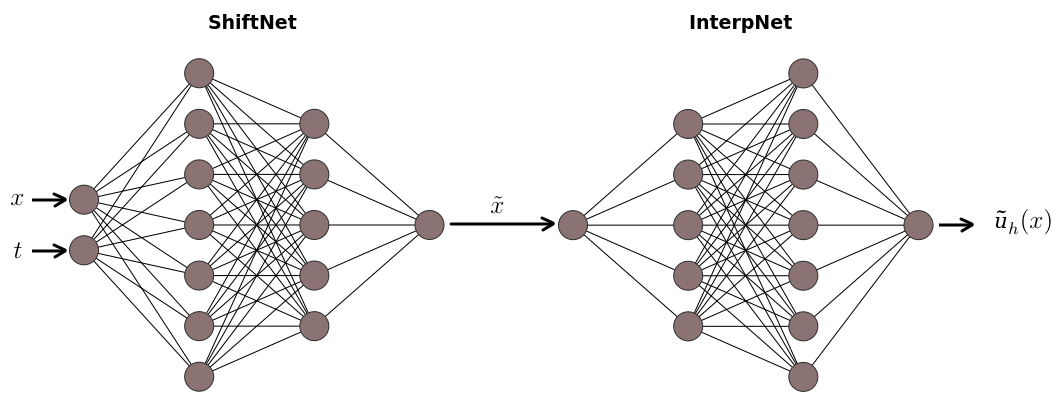}
    \caption{Deep learning architecture of NNsPOD for a $1-$dimensional model showing the continuous pipeline of information passed from ShiftNet to InterpNet.}
    \label{fig4}
\end{figure}

The introduction of a neural network with the specific task of interpolating the field values on the nearest centroids of the grid following the shift means that, once fully trained, NNsPOD is not only capable of reconstructing the values of the solution field continuously across the computational domain but, depending on how the training is performed, it will also allow for the bijective mapping to be performed for virtually unlimited non-linear advection models and in particular with those that feature high degree of variability in the \say{shape} of the snapshots collected, at full-order, in $\mathbf{M}$.
Attention must be paid in constructing such split architecture since the objective is that NNsPOD's output $\mathcal{T}$ is a sufficiently regular bijective map between the manifold and its shifted counterpart. As such, continuity in the data pipeline must be preserved from ShiftNet and InterpNet and viceversa. At the same time it is not desirable that the hyperparametric optimisation of the loss function of those networks has an effect on each other; a lack of separation of weights and biases updates lacks in fact the possibility of generalising NNsPOD transformation during online stage for those advection-dominated models that have parametric dependence in the transport field itself.
As such we devised the DNN-based algorithm, reported in \textbf{Algorithm~\ref{alg1}} to achieve such counteracting properties. The training stages of ShiftNet and InterpNet are thus separated with the latter being trained first. Once the network has learned the best possible reconstruction of the solution field of the reference configuration, its forward map will be used for the training of ShiftNet as well, in a cascaded fashion. For this reason we must optimise the loss of InterpNet (whose training set is composed by the field values of $\mathbf{u}_{\text{ref}}$) considerably more than ShiftNet's. A schematic of the structure of NNsPOD, as well as his input and output relation for the $1-$dimensional example discussed in~\ref{sub3.1} is depicted in Figure~\ref{fig4}.
We highlight the simple scalability of NNsPOD architecture and algorithm to hyperbolic equations with higher dimensional space domain. As a matter of fact ShiftNet's training set consists on the data-point's spatial distribution meaning that its input layer will have $d+P$ neurons while the output layer, feeding information to InterpNet, has $d$ neurons. 

InterpNet on the other hand will have $d$ neurons in the input layer while the output layer will always feature a single neuron, for any spatial dimension, if the numerical solution is scalar or $3$ neurons if the solution is a vector field. We finally emphasise that a discussion about the fine-tuning of the networks is out of the scope of the present work. We will nonetheless specify the full details of the architecture deployed for the shift-detection of the numerical models that follow.

\subsection{Reduction of non-uniform, non-constant linear advection equation}\label{3.3}
\begin{figure}[H]
    \centering
    \includegraphics[keepaspectratio,width=\textwidth]{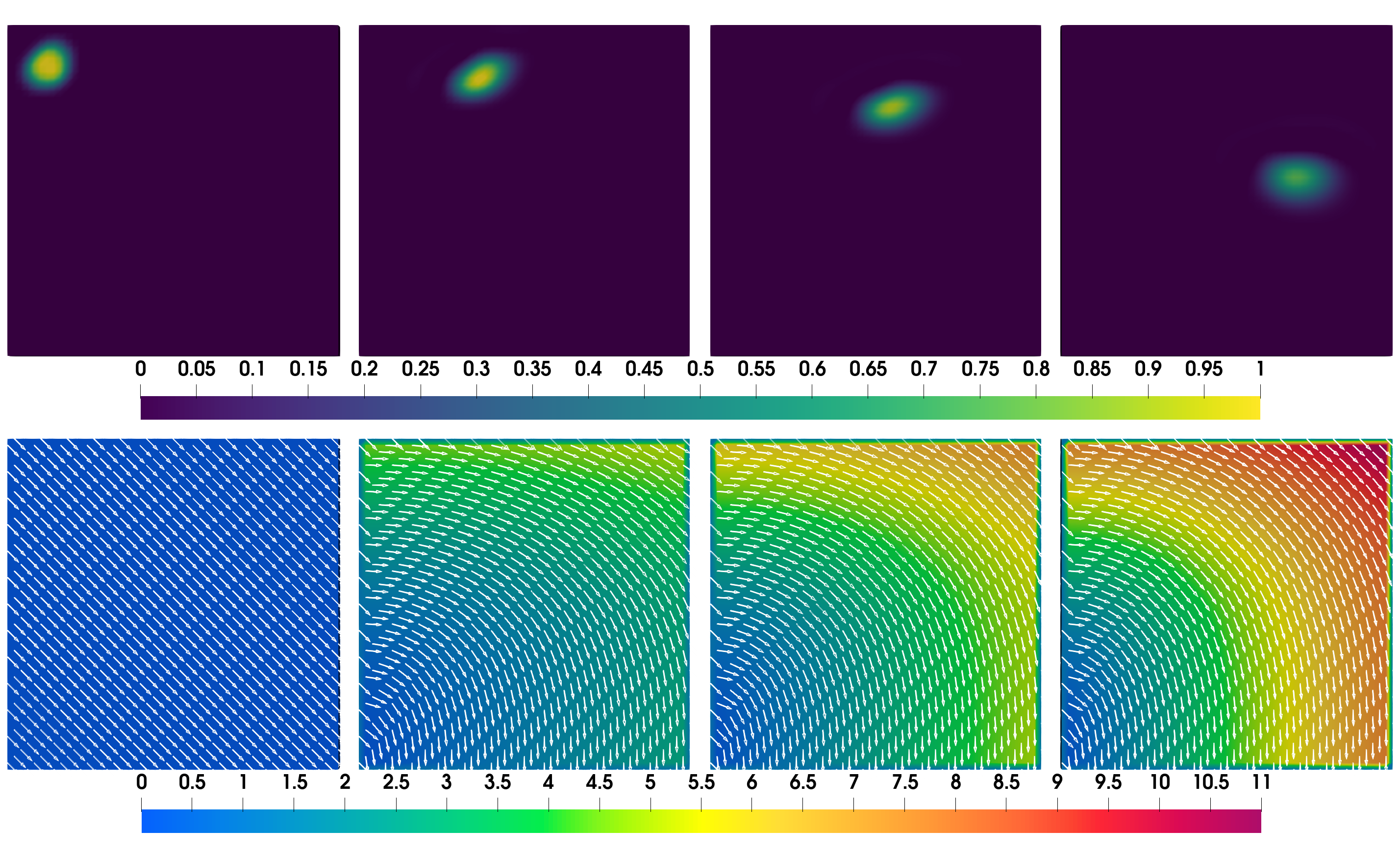}
    \caption{Different snapshots of the FOM solution (upper row) of a $2-$dimensional advection equation with non-uniform and non-constant linear transport field (bottom row).}
    \label{fig5}
\end{figure}

In order to benchmark the capability of NNsPOD automatic shift-detection the reduction of the same advection equation in $(\ref{eq3})$ is hereby derived using its workflow. A non-uniform and non-constant advection field will be used to validate the performance of the proposed algorithm; such model cannot be pre-processed by traditional shift-based methods if prior knowledge of the equation is used. As outlined in~\ref{sub2.4} in fact, a singular values spectrum analysis on sampled constant and uniform velocities will lead to an excessive computational time during the offline phase and it will also introduce a linear approximation of the characteristic curves of the hyperbolic equation. We therefore refer to the same $2-$dimensional, single parameter setting in $(\ref{eq3})$ but with $\mathbf{b}=(\frac{1}{2}\,y^2t,-2x\,t^2)$ instead; the $N_s=100$ collected snapshots are then pre-processed according to \textbf{Algorithm~\ref{alg1}}. A substantial variance in shape of the snapshots is recorded, as depicted in Figure~\ref{fig5}. this model reduction has the aim of:
\begin{itemize}
    \item Showcasing the accurate shift-based pre-processing of linear subspace approximation for the manifold of a relatively complex advection field through automatic detection i.e. using no prior knowledge of the FOM.
    \item Testing the ability of reconstructing highly \say{diffusive} snapshots collected from an hyperbolic FOM to be later applied to the more complex multiphase model in Section~\ref{Multiphase}.
\end{itemize}
The latter is a desirable quality for our data-driven algorithm to feature; highly diffusive snapshots may arise in fact not only in models with non-uniform and non-constant advection but also whenever the discretisation scheme of the convective term is of lower order of accuracy.
The settings of NNsPOD's neural networks for the shift-detection are the following
\begin{table}[H]
\centering
\begin{tabular}{|c||c|c|}
 \hline
 \multicolumn{3}{|c|}{\textbf{NNsPOD settings}} \\
 \hline
                    & \textbf{InterpNet} & \textbf{ShiftNet} \\
 \hline
\textbf{Hidden layers}$\times$\textbf{neurons} & $2\times40$ & $3\times20$ \\
\textbf{Activation function} & Sigmoid & PReLU \\
\textbf{Learning rate} & $10^{-3}$ & $10^{-4}$ \\
\textbf{Accuracy threshold} & $10^{-7}$ & $10^{-1}$ \\
 \hline
 \end{tabular}
 \end{table}

\begin{figure}[H]
     \centering
     \begin{minipage}{\textwidth}
         \begin{subfigure}[b]{0.225\textwidth}
         \centering
         \includegraphics[keepaspectratio,width=\textwidth]{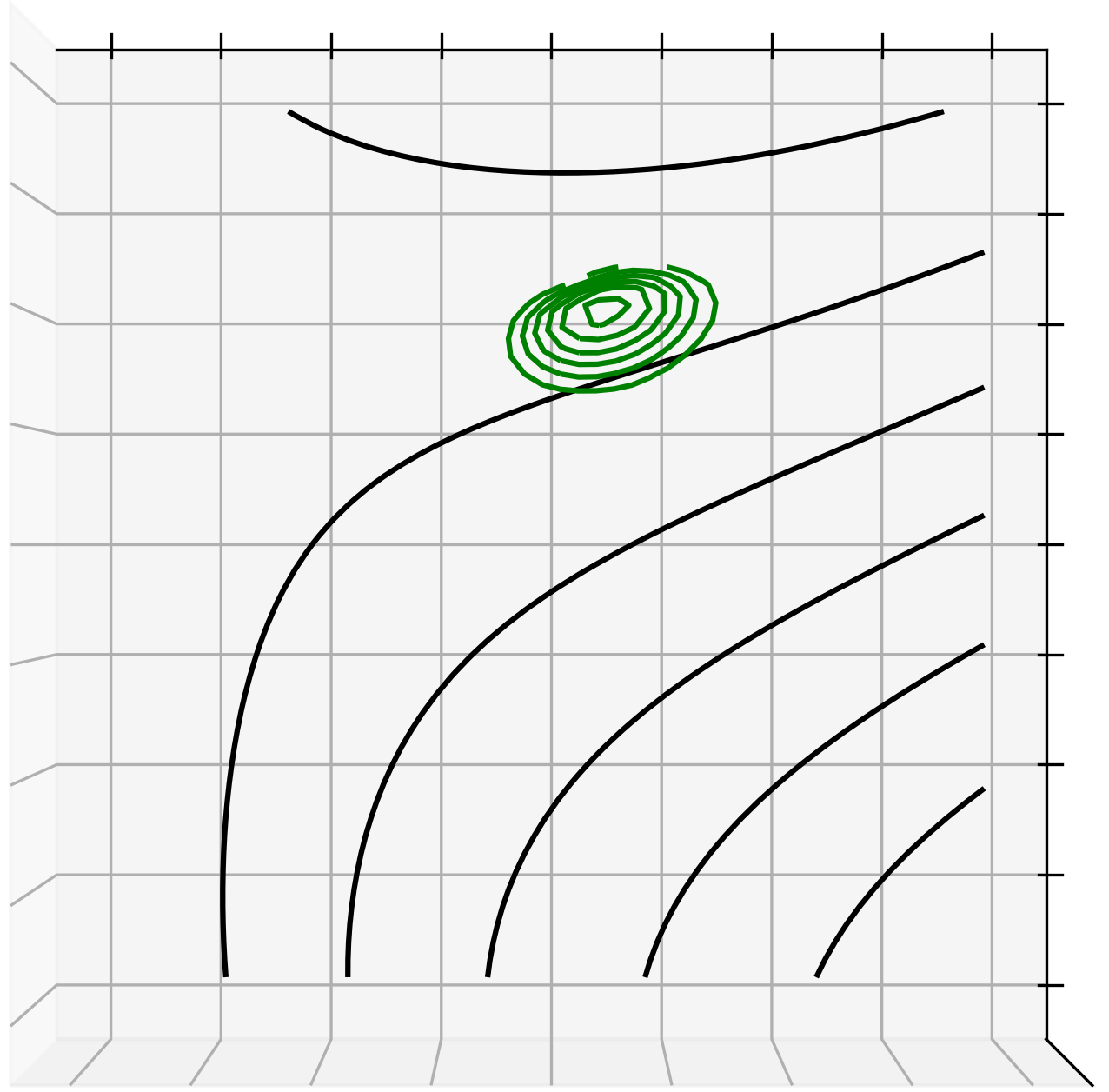}
        \end{subfigure}
        \hfill
        \begin{subfigure}[b]{0.225\textwidth}
         \centering
         \includegraphics[keepaspectratio,width=\textwidth]{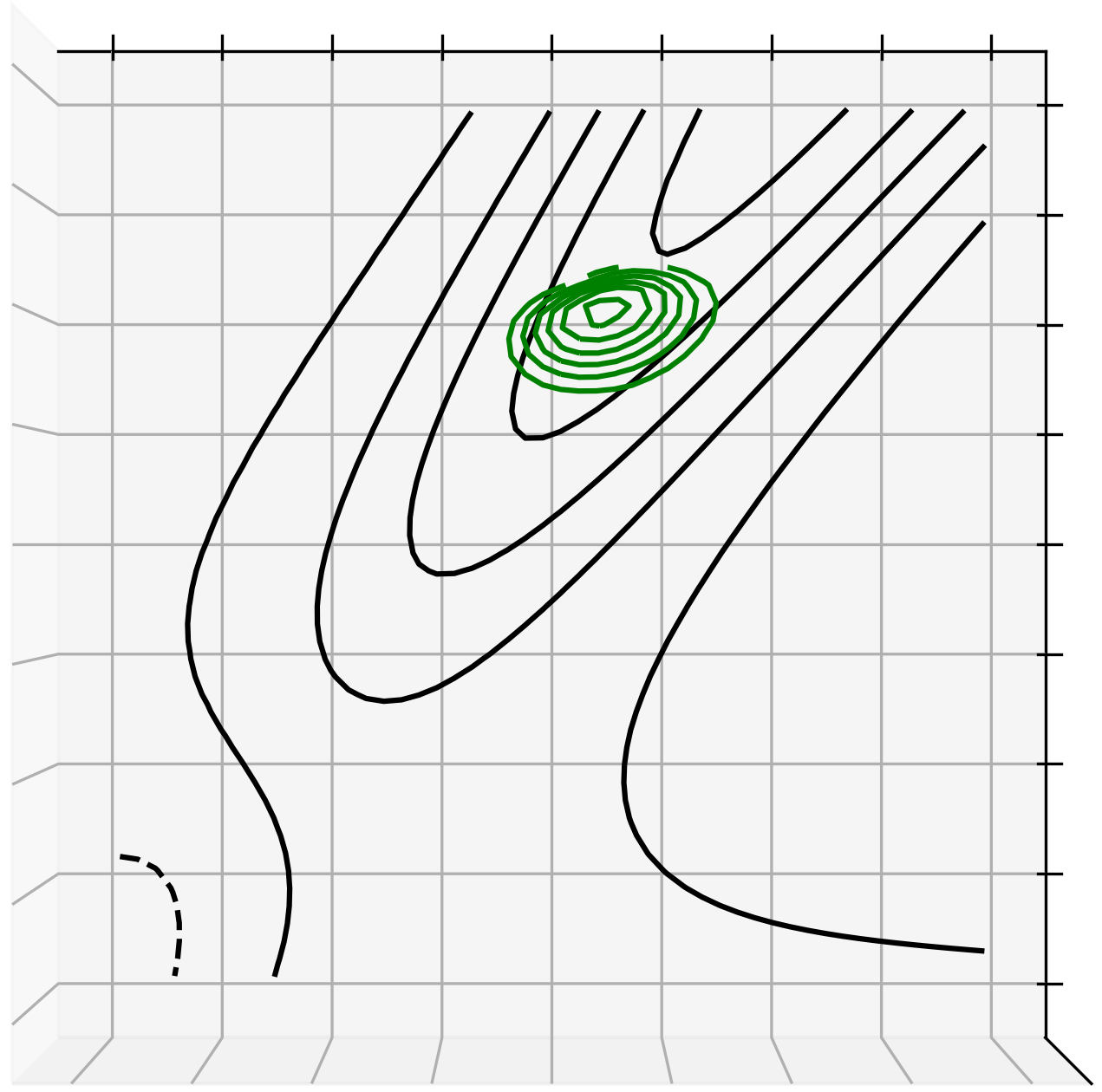}
        \end{subfigure}
        \hfill
        \begin{subfigure}[b]{0.225\textwidth}
         \centering
         \includegraphics[keepaspectratio,width=\textwidth]{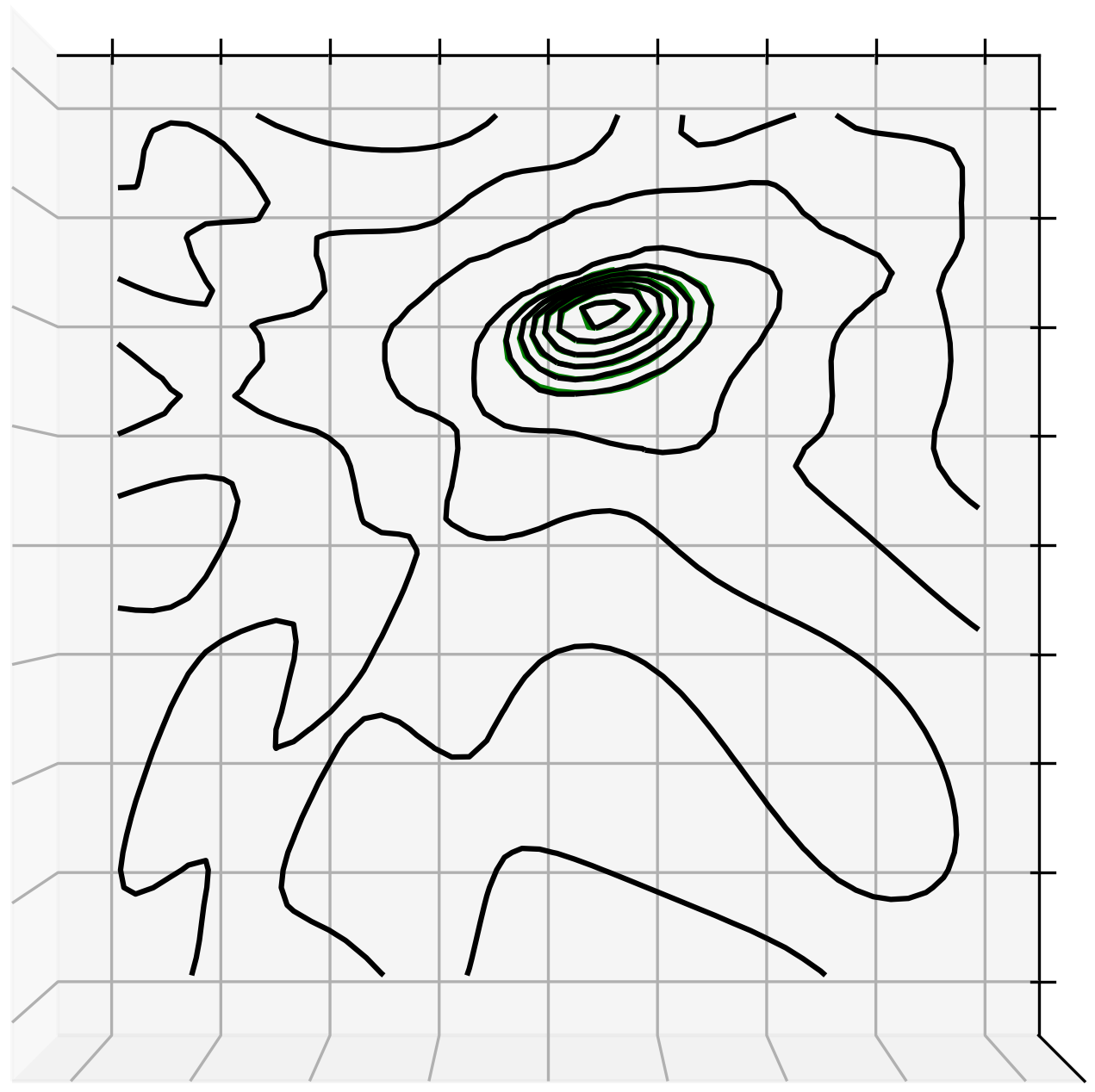}
        \end{subfigure}
        \hfill
        \begin{subfigure}[b]{0.225\textwidth}
         \centering
         \includegraphics[keepaspectratio,width=\textwidth]{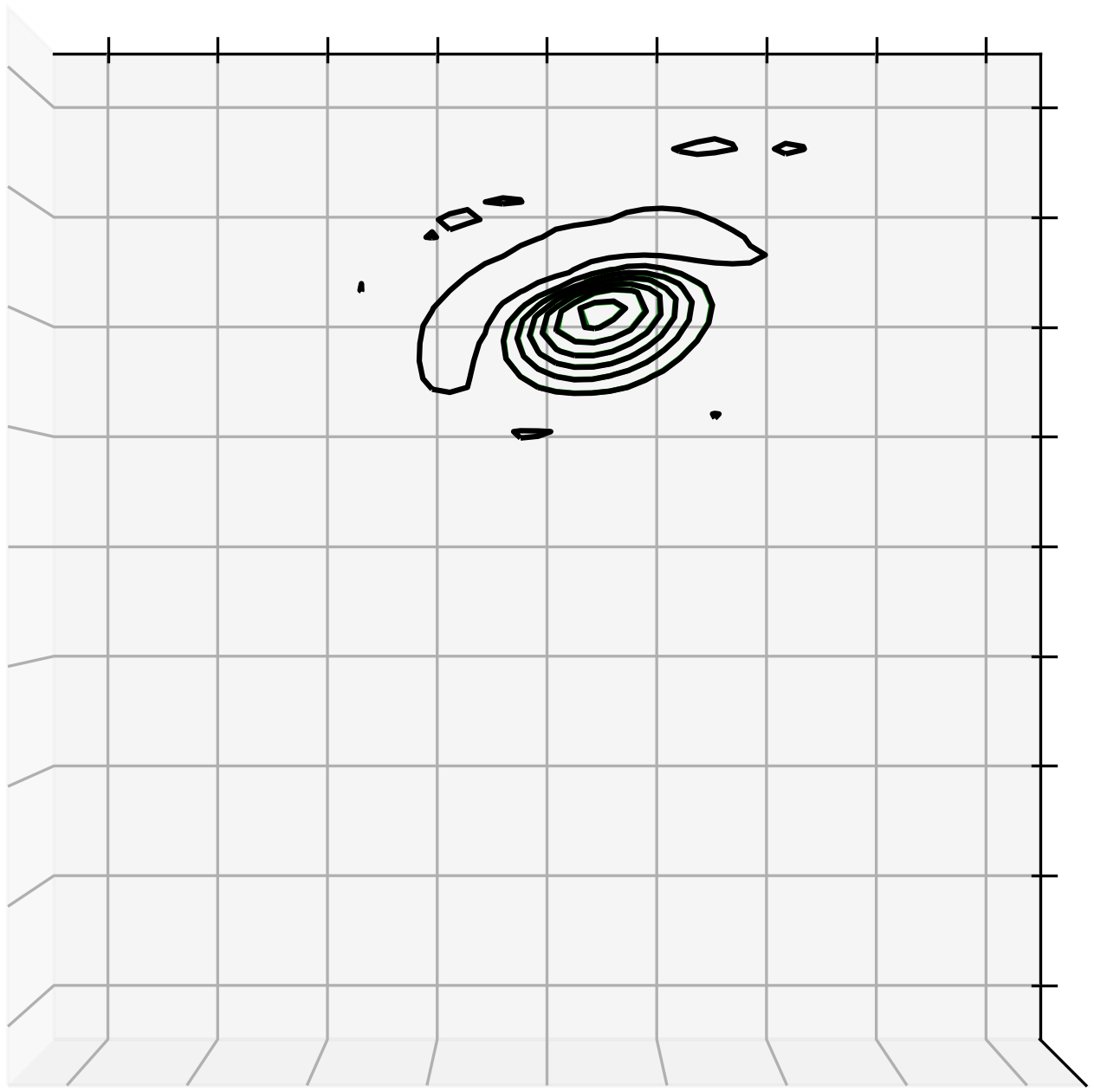}
        \end{subfigure}
     \end{minipage}
     \newline
     \newline
     \begin{minipage}{\textwidth}
         \begin{subfigure}[b]{0.225\textwidth}
         \centering
         \includegraphics[keepaspectratio,width=\textwidth]{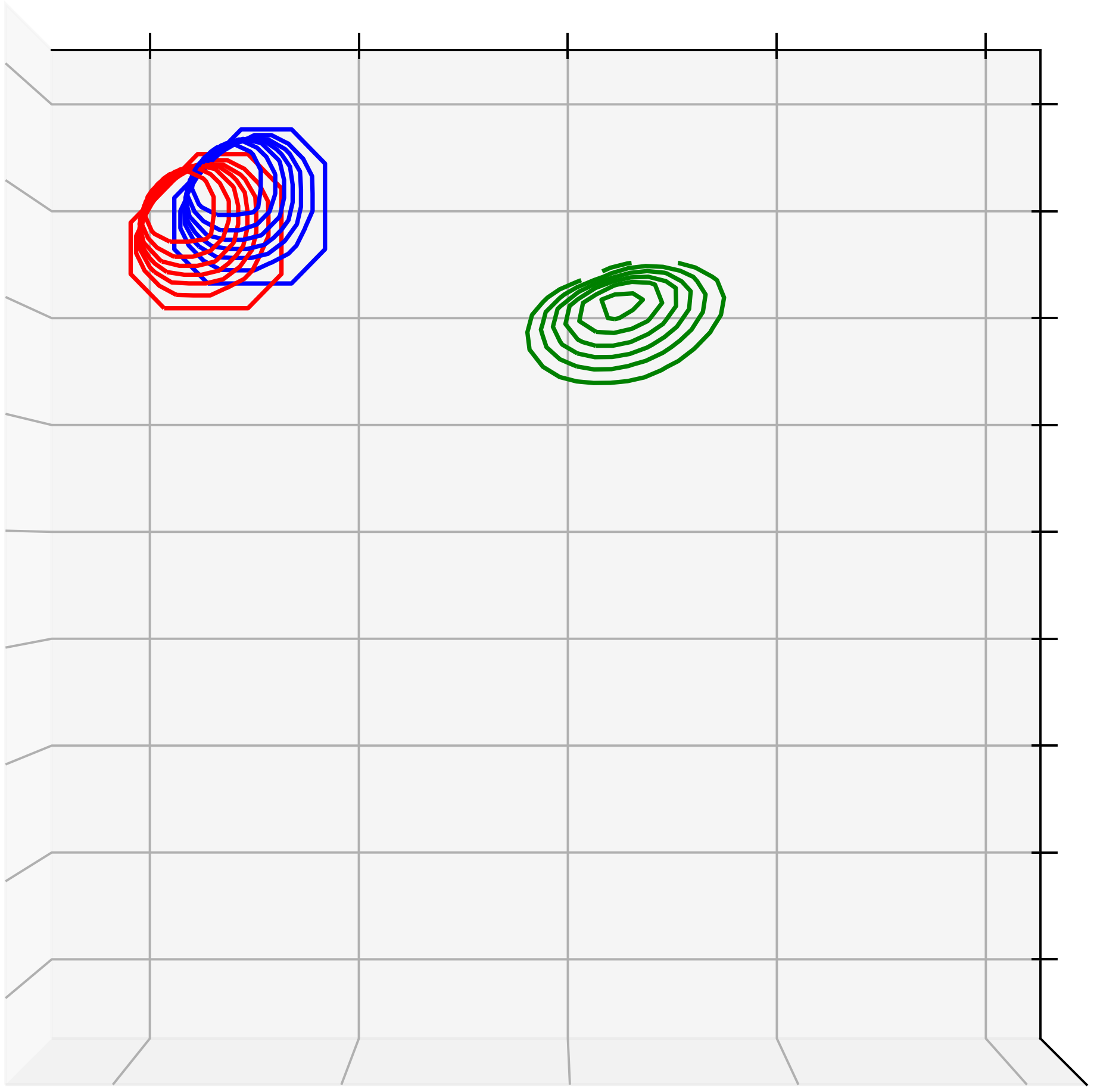}
        \end{subfigure}
        \hfill
        \begin{subfigure}[b]{0.225\textwidth}
         \centering
         \includegraphics[keepaspectratio,width=\textwidth]{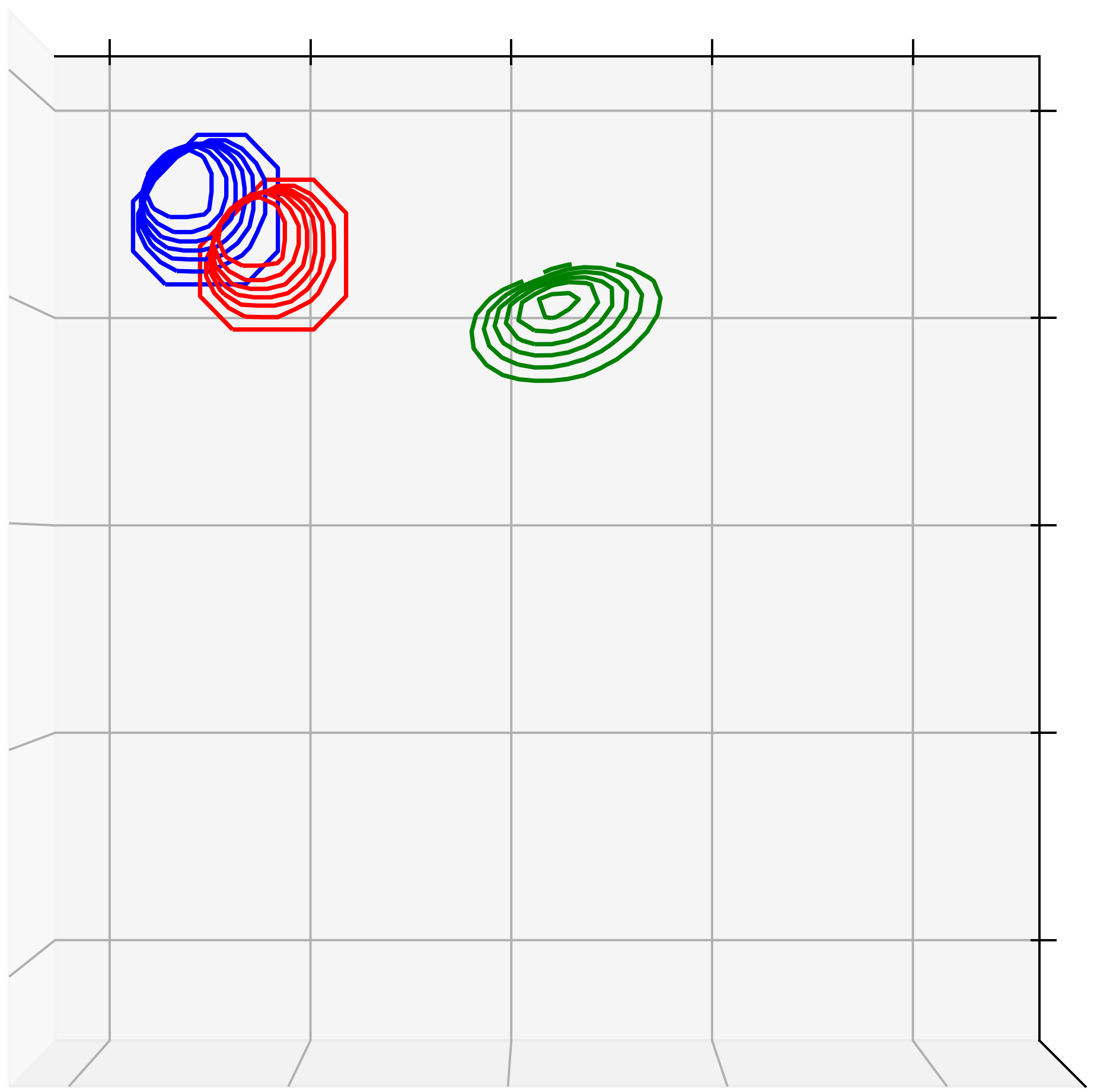}
        \end{subfigure}
        \hfill
        \begin{subfigure}[b]{0.225\textwidth}
         \centering
         \includegraphics[keepaspectratio,width=\textwidth]{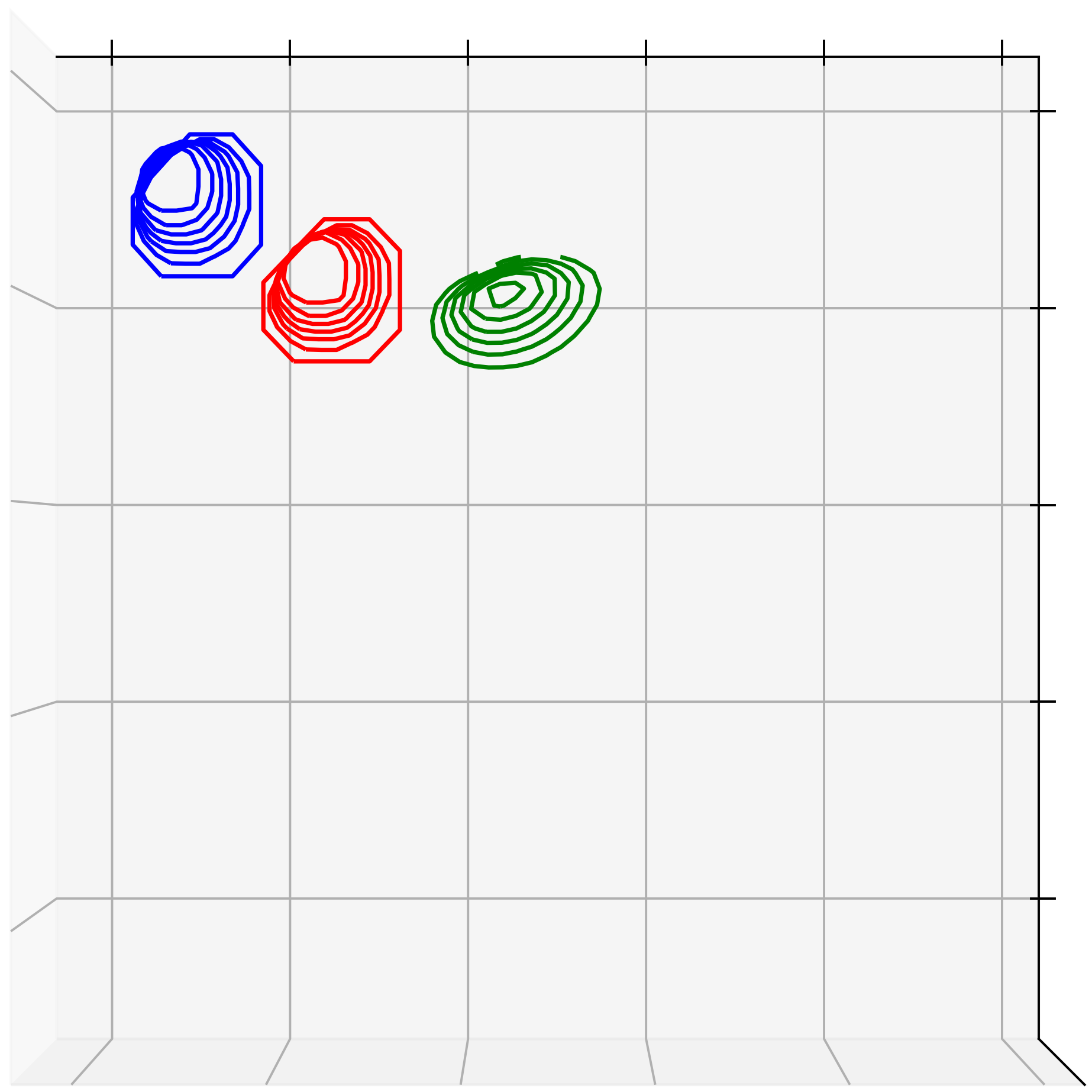}
        \end{subfigure}
        \hfill
        \begin{subfigure}[b]{0.225\textwidth}
         \centering
         \includegraphics[keepaspectratio,width=\textwidth]{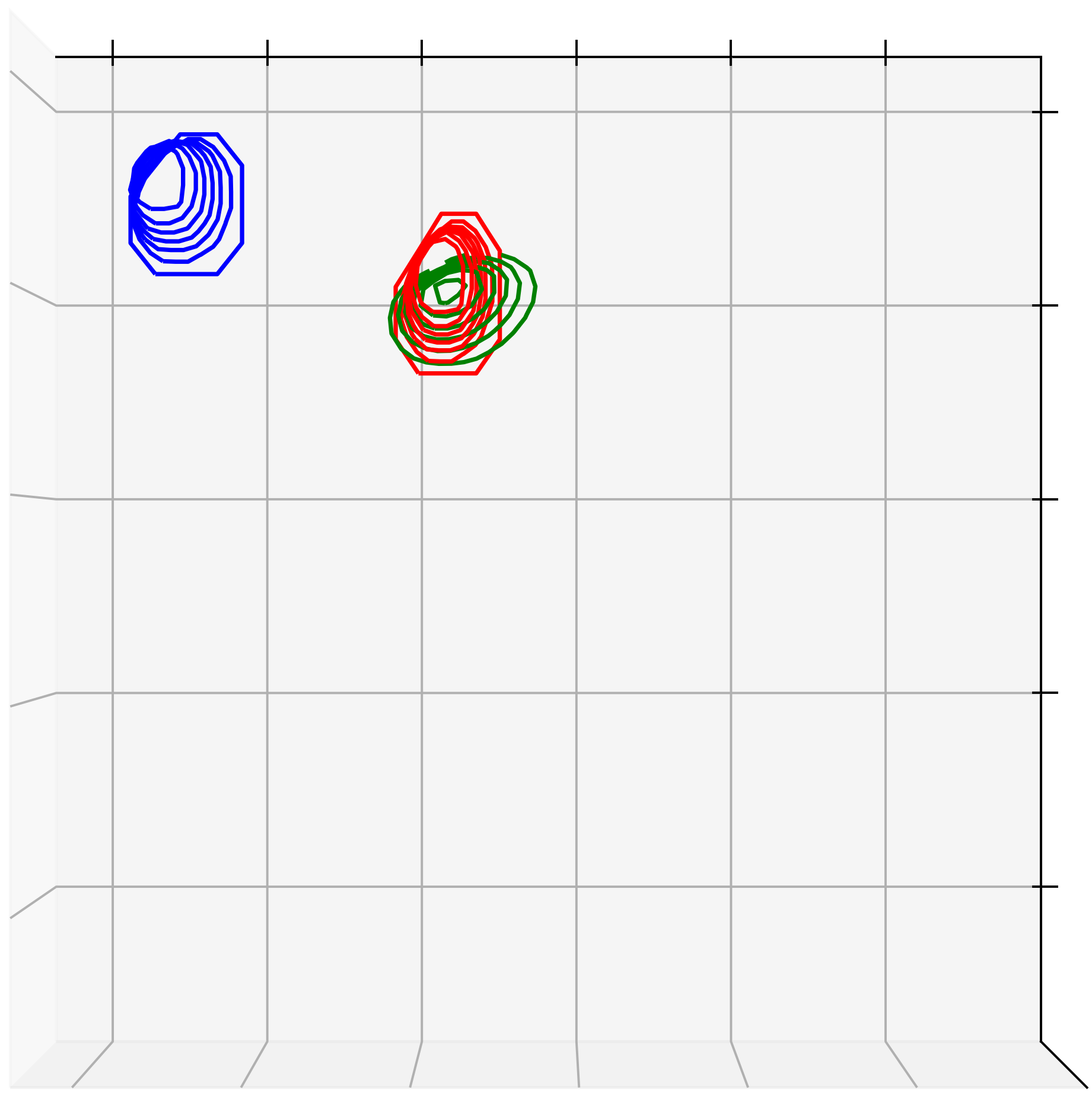}
        \end{subfigure}
     \end{minipage}
     \caption{Output of NNsPOD's split neural networks at different epochs during their separate training stages: in the upper row IntepNet's output (black contour lines) convergence to the reference configuration (green contour lines) is shown; in the bottom row the output of ShiftNet's detected coordinates (red) for a so-called test snapshot in $\mathbf{X}$ (blu) is depicted to converge to the reference configuration (green).}
     \label{fig6}
\end{figure}

The training of NNsPOD's ShiftNet neural network lasted $6$ hours on a $\text{i}7-7700$ quad processor with $16$ GB of RAM; this constituted approximately $85\%$ of the computational time of the overall offline stage (FOM snapshot collection, shift pre-processing and SVD of the shifted snapshot matrix). The different stages of the training phase of InterpNet and ShiftNet are reported in Figure~\ref{fig6}.
We tested NNsPOD performance by selecting the $80-$th snapshot as reference configuration $(\mathbf{u}_{\text{ref}}:=\mathbf{u}_h(t_{80})\in\mathbf{M})$ to highlight the possibility of choosing appropriate frames of reference for the shift-based pre-processing of $\mathcal{M}_h$, aside for the IC's one. The optimal transformation for the snapshots in $\mathbf{X}$ is thus retrieved, via automatic shift-detection, for such reference configuration.
\newpage

During the training stage we observed that many snapshots, being different in shape among each other and specifically w.r.t. the reference configuration, were stretched differently along different directions in the spatial domain by ShiftNet in order to overlap InterpNet's transformation. Nevertheless a pre-processing bijective map is derived by NNsPOD that uses no prior information regarding the mathematical model in $(\ref{eq3})$ and its physical properties (i.e. characteristic velocities).

\begin{figure}[H]
\centering
\begin{tikzpicture}

\definecolor{color2}{rgb}{0,0.75,0}
\definecolor{color1}{rgb}{0.75,0,0.75}
\definecolor{color0}{rgb}{0,0.75,0.75}

\begin{axis}[
axis on top,
legend cell align={left},
legend style={font=\tiny},
tick pos=both,
xtick style={color=black},
ytick style={color=black},
xlabel={Mode},
ymode=log,
grid=both,
legend pos=north east,
grid style={line width=.1pt, draw=gray!10},
ylabel shift = -5pt,
ylabel={Singular value (log-scale)}
]
\addplot table [color=color0, x expr=\coordindex+1, y index=0] {data/NonUnifSquare_SV_POD};
\addlegendentry{POD}
\addplot table [color=color1, x expr=\coordindex+1, y index=0] {data/NonUnifSquare_SV_sPOD};
\addlegendentry{Manually shifted POD}
\addplot table [color=color2, x expr=\coordindex+1, y index=0] {data/NonUnifSquare_SV_NNsPOD};
\addlegendentry{Automatic shift-detected POD}
\end{axis}

\end{tikzpicture}
\caption{Accuracy comparison between traditional POD and NNsPOD for a non-uniform, non-constant, linear advection field in $(\ref{eq3})$. The decay pattern follows that of an exact shift operator computed piecewise for every FOM snapshot proving the arbitrary accuracy achievable by the automatic shift-detecting bijection sought by NNsPOD.}
\label{fig7}
\end{figure}
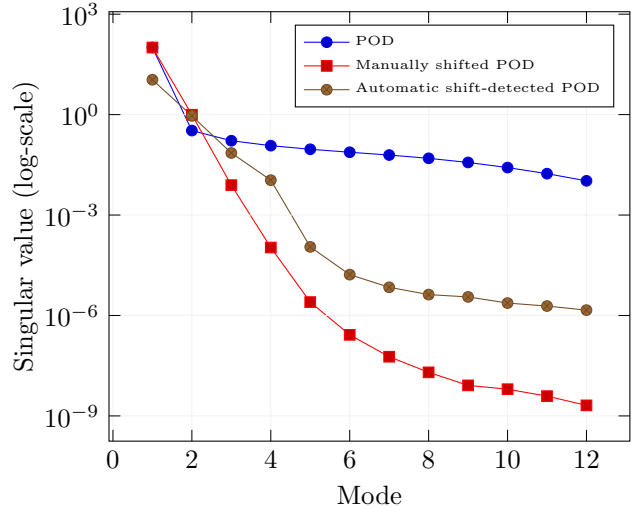

\begin{figure}[H]
\includegraphics[width=.48\textwidth]{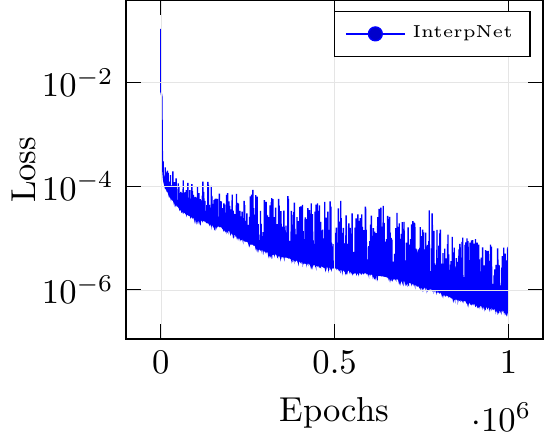}
\includegraphics[width=.48\textwidth]{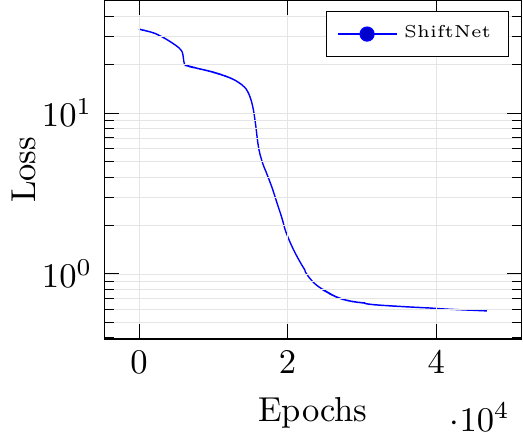}

\caption{Loss optimisation w.r.t. the training of NNsPOD's split neural networks for the IBVP $(\ref{eq3})$ with advection field $\mathbf{b} = (\frac{1}{2}y^2t,\,-2xt^2)$.}
\label{fig8}
\end{figure}

Furthermore the non-uniform and non-constant advection field has not been linearised in its parametric dependence (in this instance only w.r.t. to the time variable $t$) in order to sample the velocity space.
The SVD of the reconstructed snapshots, following NNsPOD automatically detected shift-based mapping, demonstrates a faster Kolmogorov width decay (Figure~\ref{fig7}) compared to a traditional POD algorithm; its accuracy does indeed follow the one that it would be obtained by constructing manually an exact piecewise shift operator for each of the datapoints in $\mathbf{M}$. To that regard NNsPOD can be interpreted as the statistical learning algorithm that constructs the best possible approximation of a shift operator via automatic, data-driven detection; as a result its precision is arbitrarily set by choosing the loss threshold for it to generate the bijective mapping and therefore $\tilde{\mathcal{M}_h}$. In Figure $\ref{fig8}$ the loss optimisation of InterpNet and ShiftNet are depicted for the IBVP in $(\ref{eq3})$.

\section{Reduction of a multiphase model with manifold transformation}\label{Multiphase}
We now present a numerical experiment to validate the fundamental property of
NNsPOD being applicable to non-linear advection fields. The transport of a passive scalar field in a multiphase flows of two fluids with uniform densities $\rho_1,\,\rho_2\in\mathbb{R}\,,\; \rho_1\neq\rho_2$ is chosen as FOM to generalise the reduction of a $2-$dimensional linear hyperbolic equation as addressed in Section \ref{NNsPOD}. 

The reason motivating this choice lies within the mathematical model
describing the multiphase flow itself being the coupling between the (incompressible) Navier-Stokes momentum balance equation and an advection equation. The coupling stems from the velocity field, which is unknown and derived numerically as a solution of the Navier-Stokes equation; the divergence of such field ultimately becomes the transport operator for a scalar field $\alpha(\mathbf{x},t,\boldsymbol{\mu})$ that models the pointwise fraction of volume of the two fluids within the cells of the discretised computational domain.
Furthermore, being Navier-Stokes a non-linear PDE, the no-prior knowledge quality of NNsPOD is thereby assessed as its ability of deducing the proper pre-processed linear manifold approximation with a generalisation to non-linear advection fields.

\subsection{The Volume-of-Fluid method for Eulerian interface tracking}
A multiphase flow is modelled by two fluids occupying a certain fraction of the volume domain. The time evolution of each fluid is governed by the unsteady Navier-Stokes equations
\begin{equation}\label{eq8}
    \begin{cases}
        \partial_t\,\rho + \nabla\cdot(\rho\mathbf{U}) = 0\,,\\
        \partial_t(\rho\mathbf{U}) + \nabla\cdot(\rho\mathbf{U}\otimes\mathbf{U}) = -\nabla p + \nabla\cdot\boldsymbol{\tau} + \mathbf{F}\,,
    \end{cases}
\end{equation}
in which the first equation describes the mass conservation principle (continuity equation) whereas the second models Netwon's second law for the conservation of (linear) momentum. Assuming incompressible regime and absence of external forces $(\mathbf{F} = \mathbf{0})$ acting of the fluids $(\ref{eq8})$ reduces to
\begin{equation}\label{eq9}
    \begin{cases}
        \nabla\cdot\mathbf{U} = 0\,,\\
        \partial_t\mathbf{U} + (\mathbf{U}\cdot\nabla)\mathbf{U} = -\frac{1}{\rho}\nabla\,p + \nu\Delta\,\mathbf{U}\,.
    \end{cases}
\end{equation}
The system above is coupled in the pressure $p(\mathbf{x},t)$ and velocity $\mathbf{U}(\mathbf{x},t)$ fields which are the unknowns of the problem at hand that has to be solved numerically. We observe that $(\ref{eq9})$ describes the dynamics of one single fluid with (uniform) volumetric density $\rho$ and kinematic viscosity $\nu$. In the multiphase flow however we have two fluids, separated by a sharp interface, that evolve in a coupled fashion with each other.
  By assuming that the two fluids are immiscible we can describe the physical properties of the flow as averaged across the domain; to that end the fraction of volume fields $\alpha_1,\,\alpha_2\in[0,1]$ are introduced which are defined through a constitutive relation $\alpha_1\Omega + \alpha_2\Omega = \Omega \Rightarrow \alpha_1 + \alpha_2 = 1 \Rightarrow \alpha := \alpha_1 = 1 - \alpha_2$ that acts as a constraint on the model. The new density field is thus defined as $\alpha(\mathbf{x},t) = \alpha(\mathbf{x},t)\rho_1 + (1-\alpha(\mathbf{x},t))\rho_2$ and its time evolution is governed by an advection equation coupled with the constitutive relation derived above
\begin{equation}\label{eq10}
    \begin{cases}
        \alpha\rho_1 + (1-\alpha)\rho_2 = \alpha \;,\;\;\alpha\in[0,1] \;\:\forall\mathbf{x}\in\Omega \,,\\
        \partial_t\,\alpha + \nabla\cdot(\mathbf{U}\alpha) = 0\,,
    \end{cases}
\end{equation}
with $\mathbf{U}(\mathbf{x},t)$ being the (non-linear) velocity field in $(\ref{eq9})$ which is thus coupled with $(\ref{eq10})$. Eulerian multiphase modelling tracks the instantaneous and pointwise changes in the sharp interface between the two fluids; continuous-continuous and dispersed-continuous phase interaction are the two large families of algorithms used to model multiphase flow. The latter is used in simulating dispersed particles (solid phases), droplets (liquid phases) and bubbles (gaseous phases) within a larger  continuous fluid phase; the former instead is used whenever there is a presence of a continuous sharp interface between two fluids. 
The presence of an additional transport operator in the advection equation for field $\alpha(\mathbf{x},t)$ suggests that traditional bounded schemes for its numerical discretisation might not suffice. The operator becomes in fact highly diffusive in proximity of the interface in a processed called excessive smearing. One possibility to mitigate the phenomena, as outlined in \cite{renardy_numerical_2001}, is to introduce a numerical compression term $\nabla\cdot(\alpha\,(1-\alpha)\,\mathbf{U}_r)$ in the advection equation in $(\ref{eq10})$ built ad-hoc s.t. it takes null values everywhere in the domain except in proximity of the interface where $\alpha\in(0,1)$ (being $\mathbf{U}_r :=C_{\alpha}||\mathbf{U}||\frac{\nabla\alpha}{||\nabla\alpha||}$ with $C_{\alpha}\in[0,1]$ as a free parameter). 
The resulting system from $(\ref{eq10})$ thus is
\begin{equation}\label{eq11}
    \begin{cases}
            \alpha\rho_1 + (1-\alpha)\rho_2 = \alpha \;,\;\;\alpha\in[0,1] \;\:\forall\mathbf{x}\in\Omega\,, \\
        \partial_t\,\alpha + \nabla\cdot(\alpha(\mathbf{U} + (1-\alpha)\mathbf{U}_r)) = 0\,.
        \end{cases}
\end{equation}
Then, a discretisation scheme is devised for the advection terms, with the compressive interface capturing (CICSAM) and piecewise linear interface construction (PLIC) being the most widespread (we refer to \cite{okagaki_numerical_2021} for a detailed explanation and derivation).
\subsection{Full-order multiphase numerical modelling}
\begin{figure}[H]
    \centering
    \includegraphics[keepaspectratio,width=\textwidth]{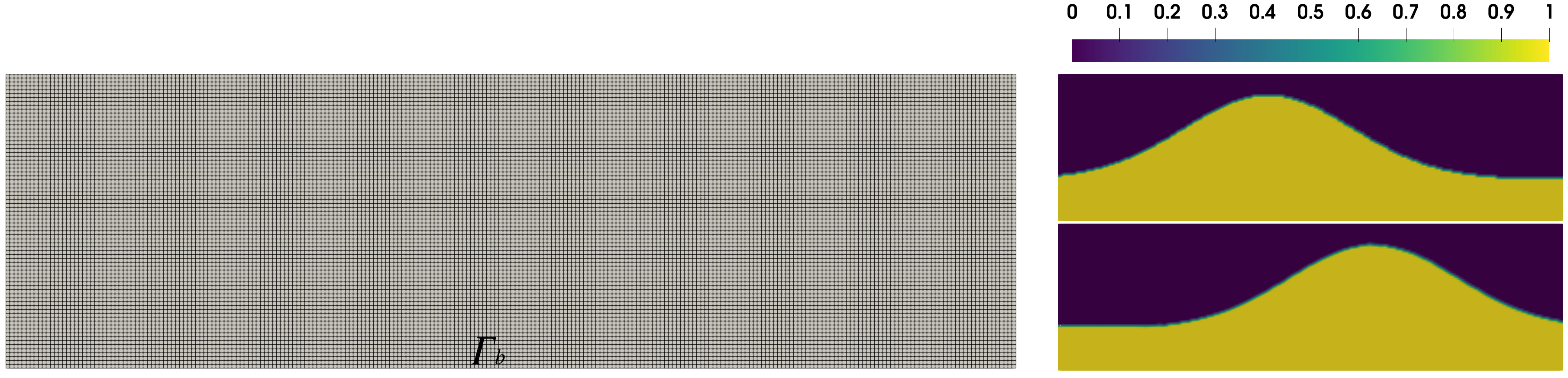}
    \caption{Illustration of the FOM for the IBVP in $(\ref{eq12})$: on the left the discrete computational grid is shown with the indication for the bottom partition of the boundary where the no-slip condition is applied; on the right two snapshots of the solution field $\boldsymbol{\alpha}_h$ are reported (IC on the top, $100-$th snapshot on the bottom).}
    \label{fig9}
\end{figure}

To benchmark the performance of NNsPOD automatic shift-detection and construction of pre-processed linear manifold approximation for the non-linear advection of an hyperbolic equation, a simplified model of a multiphase simulation is solved numerically. While retaining the parametric dependence on the density and viscosity for the couple of fluids, for the sake of simplicity we restrict the parameter space, paired to the FOM, to be the exclusively the time variable as previously derived for the linear models in Section \ref{Pre-processing} and \ref{NNsPOD}. To this end we consider the following IBVP for the time-evolution of a sharp interface separating water $(\rho_1=10^3\,,\;\nu_1=10^{-6})$ from air $(\rho_2=1\,,\;\nu_2=1.48\cdot10^{-5})$
\begin{equation}\label{eq12}
    \begin{cases}
        \nabla\cdot\mathbf{U} = 0\,, \\
        \partial_t\mathbf{U} + (\mathbf{U}\cdot\nabla)\mathbf{U} = -\frac{1}{\rho}\nabla\,p + \nu\Delta\,\mathbf{U}\,, \\
        \mathbf{U} = (0,0)\;,\;\;\forall\mathbf{x}\in\Gamma_b\,,\quad\partial_n\,\mathbf{U} = 0\;,\;\forall\mathbf{x}\in\partial\Omega\setminus\Gamma_b\,, \\
        \mathbf{U}(\mathbf{x},0) = (0.25, 0)\,,\;\;\forall\mathbf{x}\in\Omega\,, \\
        \alpha\rho_1 + (1-\alpha)\rho_2 = \alpha \;,\;\;\alpha\in[0,1] \;\:\forall\mathbf{x}\in\Omega \,,\\
        \partial_t\,\alpha + \nabla\cdot(\alpha(\mathbf{U} + (1-\alpha)\mathbf{U}_r)) = 0\,,\;\;C_{\alpha}=1\,, \\
        \alpha = 1\,,\;\;\forall\mathbf{x}\in\Gamma_b\,,\quad\nabla\alpha = 0\;,\quad\forall\mathbf{x}\in\partial\Omega\setminus\Gamma_b\,, \\
        \alpha(\mathbf{x}, 0) = 1 \;,\quad\forall\mathbf{x}\in\Omega\;\;\text{s.t.}\;\;y < e^{-\frac{x^2}{2}}\,,
    \end{cases}
\end{equation}
where the IC for field $\alpha(\mathbf{x},t)$, which identifies the initial configuration of interface between the two fluids, features a Gaussian profile. The BCs are set to comply with the no-slip condition for the bottom partition of the boundary $\Gamma_b$ whereas a null gradient through top, the inlet and outlet partitions is allowed (see Figure \ref{fig9} for reference). The setting is intended to simulate a multiphase flow through a $2-$dimensional canal with open lid at atmospheric pressure. The domain $\Omega = [-2.5, 3.5]\times[-0.5, 1.25]$ is discretised in a homogeneous set of collocated rectangular cells with a computational grid of $250\times75$ centroids $(N_h = 18,750)$. The time interval for the simulation is set to $t\in[0, 5]$ which is in turn discretised into $100$ steps to comply with the CFL condition thereby leading to the following manifold definition
\begin{equation*}
    \mathcal{M}_h := \big\{\boldsymbol{\alpha}_h(t_k)\in\mathcal{V}_h\,,\;k=1,\dots,100\big\}
\end{equation*}
In the statistical learning configuration of NNsPOD we remark that the offline snapshot collection generates a dataset $\mathbf{M}\in\mathbb{R}^{N_s\times N_h}$, associated to $\mathcal{M}_h$, of cardinality $N_s = 100$ and dimensionality (number of features) $N_h = 18.75\cdot10^3$.

 \begin{figure}[H]
     \centering
     \begin{minipage}{\textwidth}
         \begin{subfigure}[b]{0.225\textwidth}
         \centering
         \includegraphics[keepaspectratio,width=\textwidth]{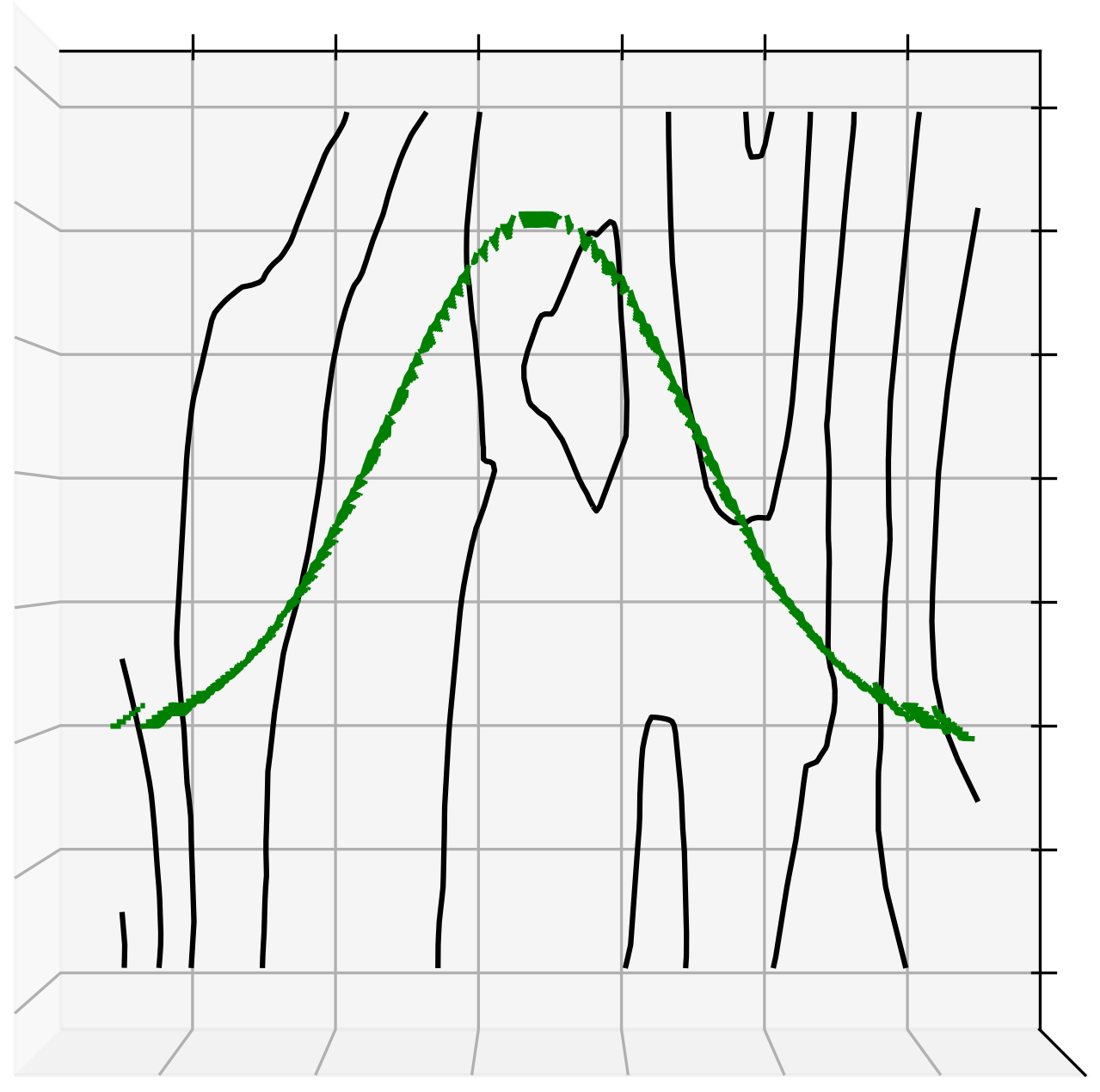}
        \end{subfigure}
        \hfill
        \begin{subfigure}[b]{0.225\textwidth}
         \centering
         \includegraphics[keepaspectratio,width=\textwidth]{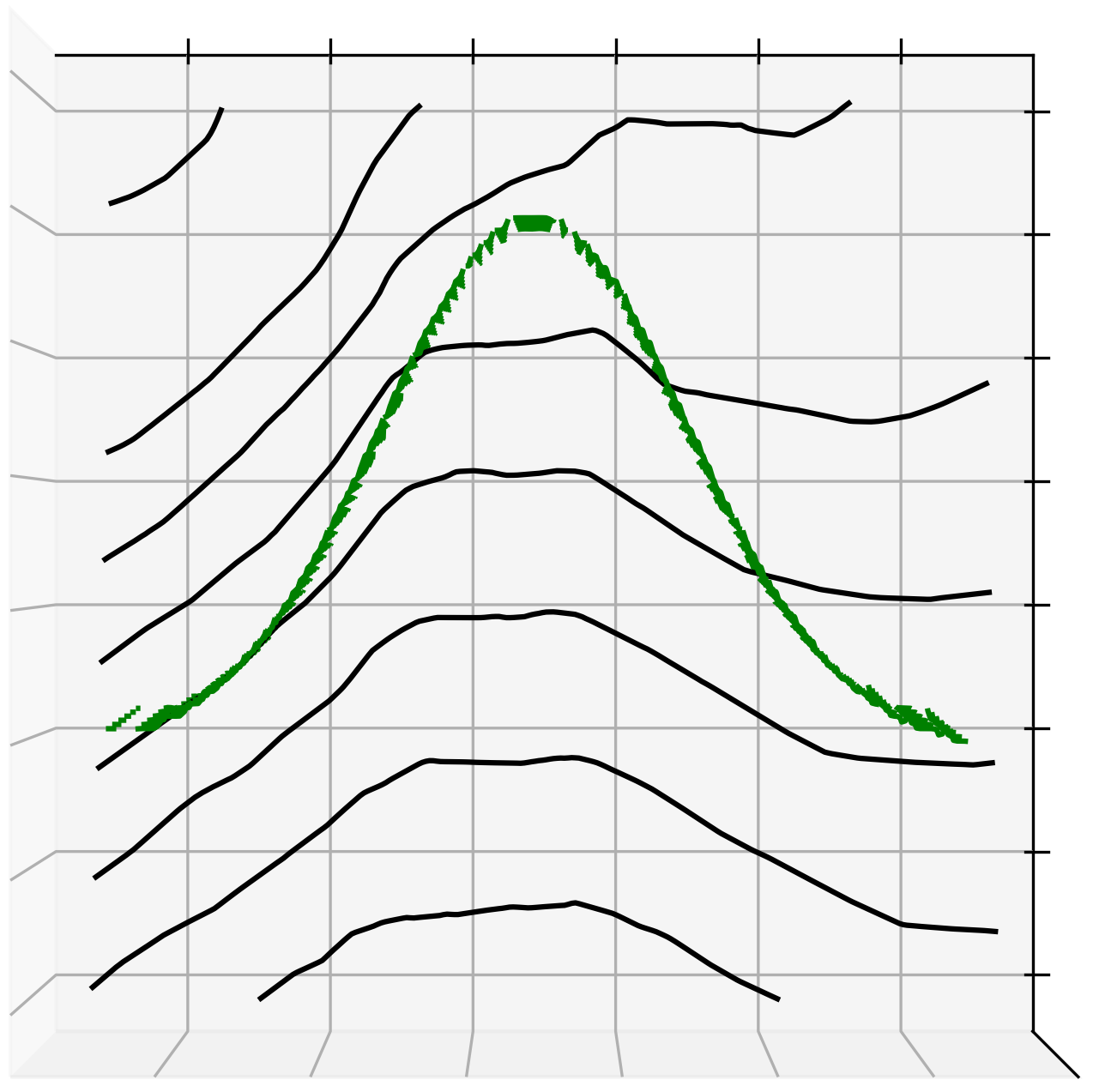}
        \end{subfigure}
        \hfill
        \begin{subfigure}[b]{0.225\textwidth}
         \centering
         \includegraphics[keepaspectratio,width=\textwidth]{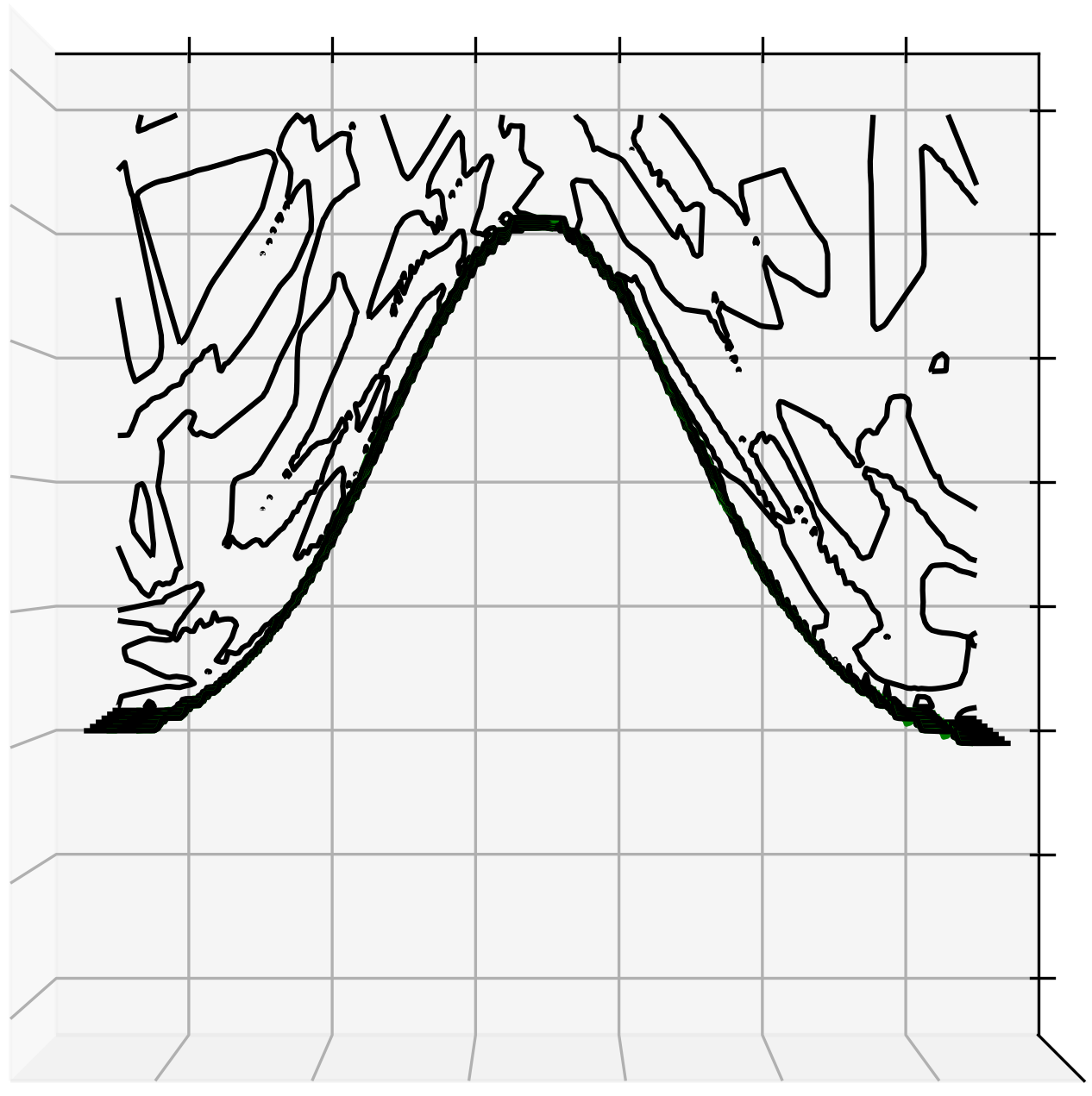}
        \end{subfigure}
        \hfill
        \begin{subfigure}[b]{0.225\textwidth}
         \centering
         \includegraphics[keepaspectratio,width=\textwidth]{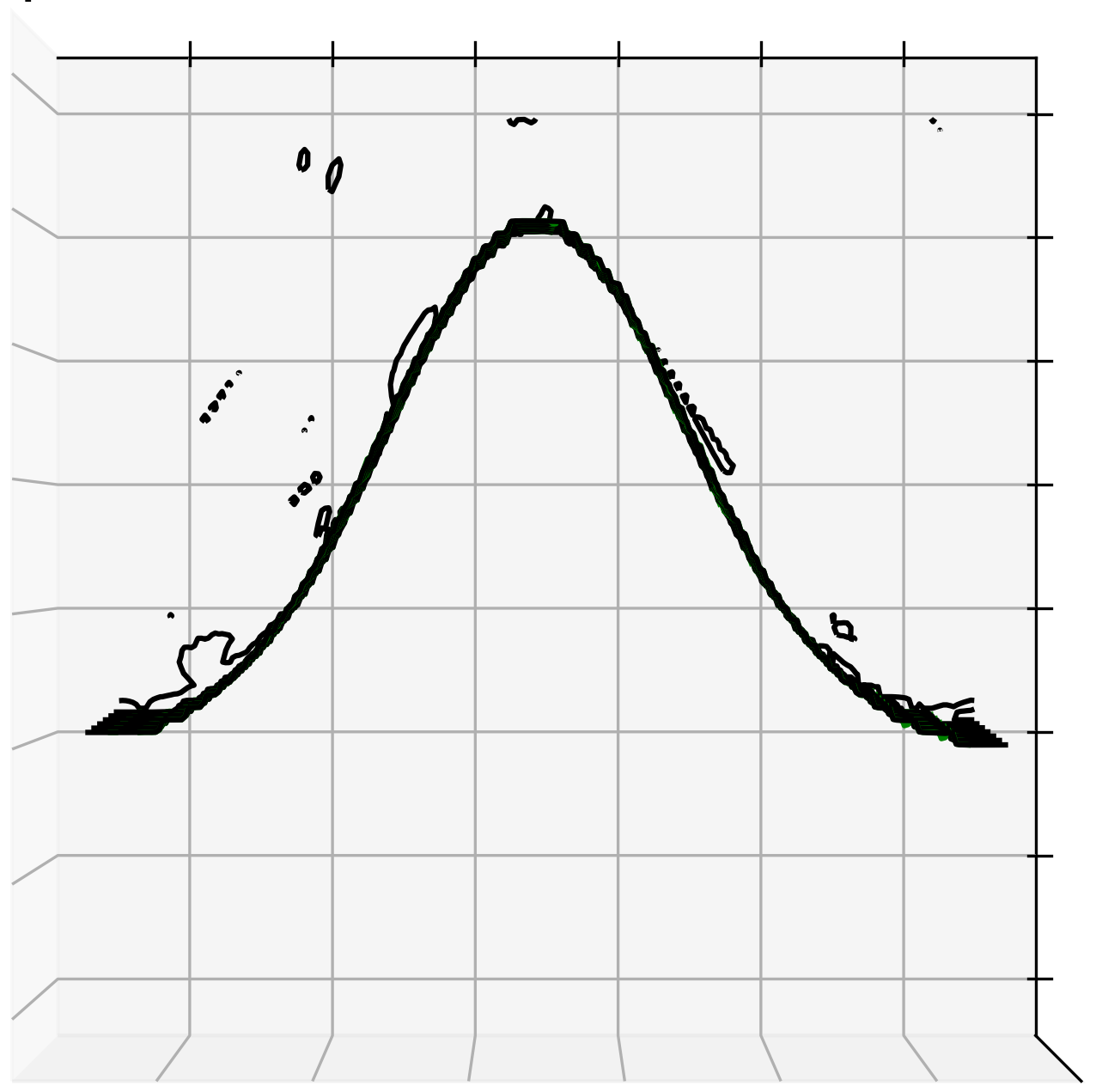}
        \end{subfigure}
     \end{minipage}
     \newline
     \newline
     \begin{minipage}{\textwidth}
         \begin{subfigure}[b]{0.225\textwidth}
         \centering
         \includegraphics[keepaspectratio,width=\textwidth]{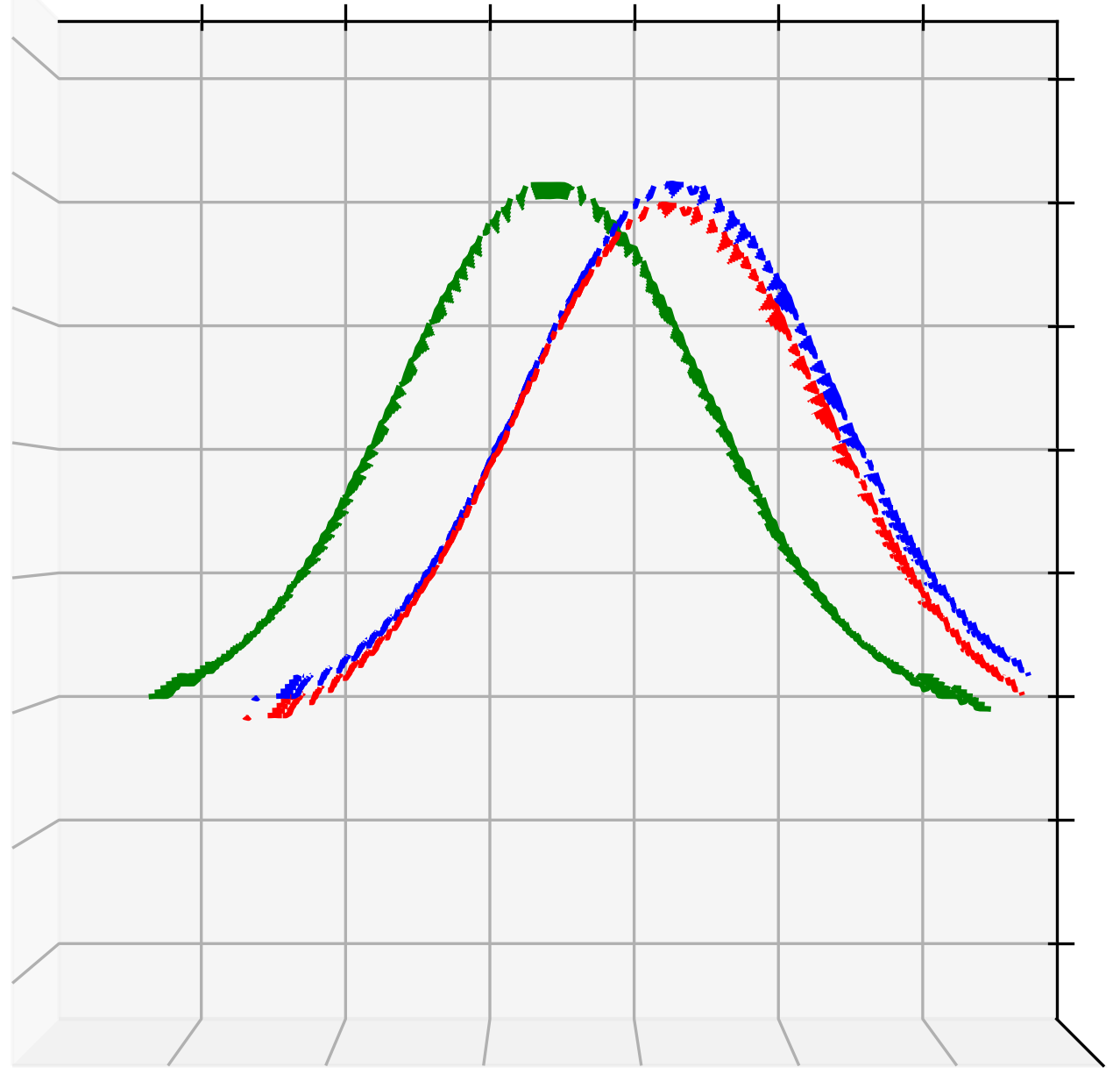}
        \end{subfigure}
        \hfill
        \begin{subfigure}[b]{0.225\textwidth}
         \centering
         \includegraphics[keepaspectratio,width=\textwidth]{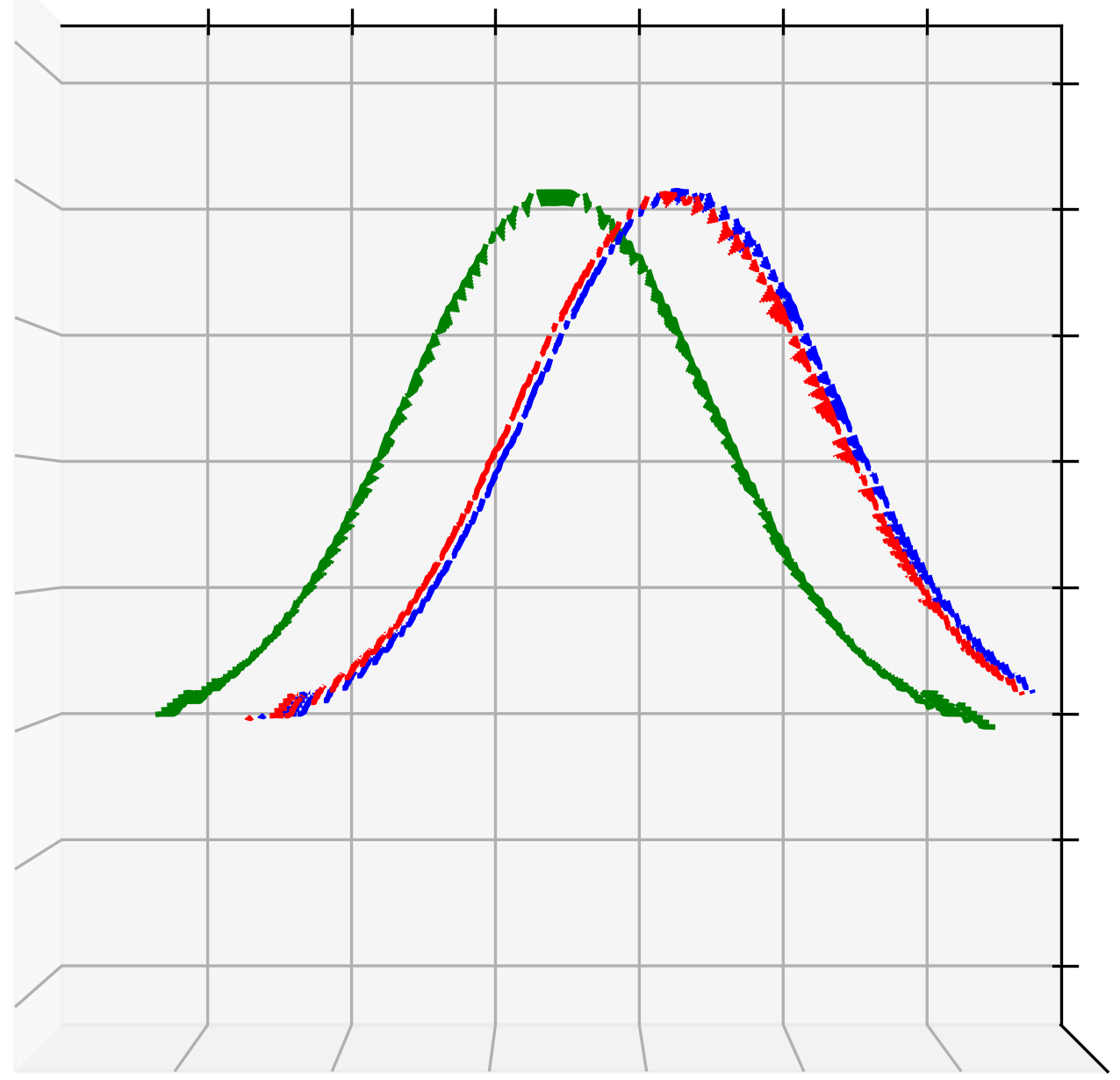}
        \end{subfigure}
        \hfill
        \begin{subfigure}[b]{0.225\textwidth}
         \centering
         \includegraphics[keepaspectratio,width=\textwidth]{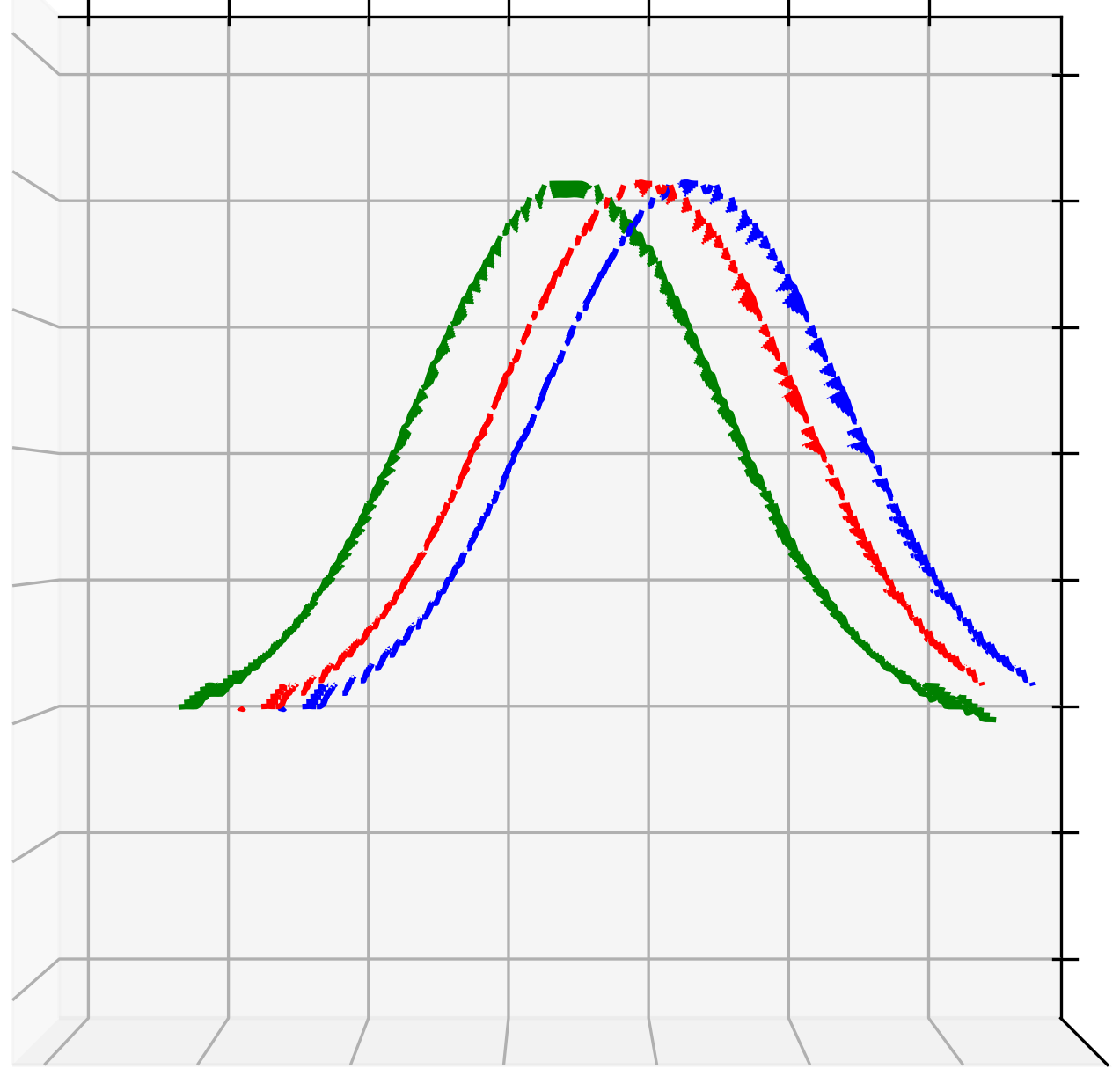}
        \end{subfigure}
        \hfill
        \begin{subfigure}[b]{0.225\textwidth}
         \centering
         \includegraphics[keepaspectratio,width=\textwidth]{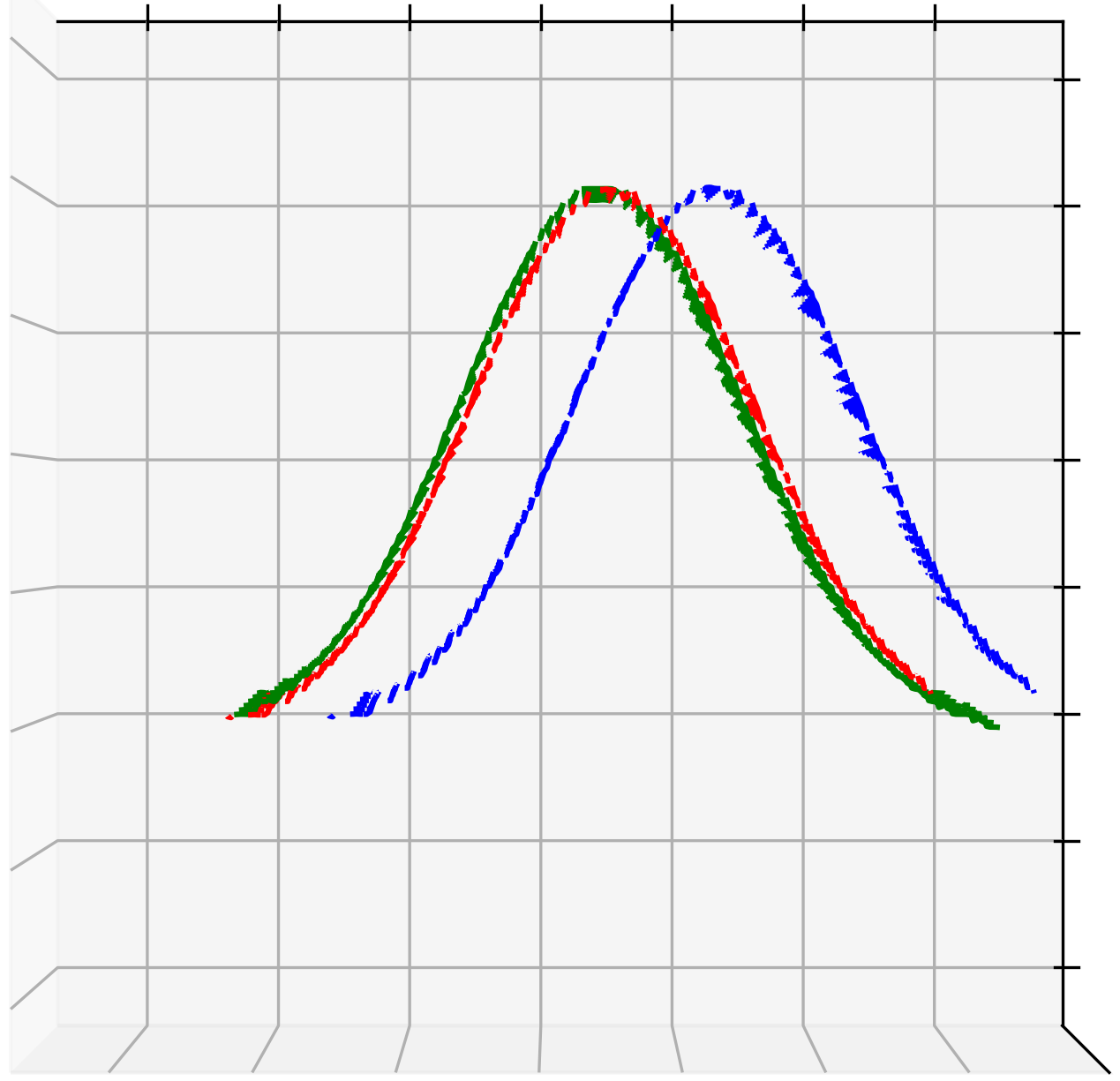}
        \end{subfigure}
     \end{minipage}
     \caption{NNsPOD's neural networks at different epochs during the training for the multiphase FOM in $(\ref{eq10})$: the upper row displays IntepNet's (black) convergence to $\boldsymbol{\alpha_{\text{ref}}}$ (green) while the bottom row shows the convergence of ShiftNet's output (red) for test snapshot in $\mathbf{X}$ (blu).}
     \label{fig10}
\end{figure}

\subsection{Automatic shift-detection and linear manifold reconstruction}
The shift-detection of the algorithm regards the correct reconstruction of the interface separating the two phases of fluids w.r.t. a reference configuration for the problem. As discussed in Section \ref{NNsPOD} the choice for an appropriate snapshot for the training of InterpNet does not affect the ability of ShiftNet in constructing a pre-processed manifold approximation therefore, even for highly diffusive transports, one can always select any configuration of the numerical field $\boldsymbol{\alpha}_h$ for the subsequent automatic shift detection. Here we selected the $35-$th snapshot to be the reference configuration $\boldsymbol{\alpha}_{\text{ref}}:=\boldsymbol{\alpha}_h(t_{35})\in\mathbf{M}$. The split architecture of NNsPOD was set with the following parameters.
\begin{table}[H]
\centering
\begin{tabular}{|c||c|c|}
 \hline
 \multicolumn{3}{|c|}{\textbf{NNsPOD settings}} \\
 \hline
                    & \textbf{InterpNet} & \textbf{ShiftNet} \\
 \hline
\textbf{Hidden layers}$\times$\textbf{neurons} & $4\times40$ & $5\times25$ \\
\textbf{Activation function} & HardSigmoid & PReLU \\
\textbf{Learning rate} & $10^{-5}$ & $10^{-6}$ \\
\textbf{Acccuracy threshold} & $10^{-4}$ & $10^{\,2}$ \\
 \hline
 \end{tabular}
 \end{table}
 The field $\boldsymbol{\alpha}_h$ targeted by the shift-detection, and therefore involved in the minimisation of the loss function, has values in $[0,1]$; intuitively in most of the cells of the FOM its values is either one of the infimum (i.e. air only) or the supremum (i.e. water only) with the exception of the interface where the gradient has non-null values. As such, in order for NNsPOD to emulate such behaviour, the hard sigmoid activation function for InterpNet's neurons, which is responsible for \textsl{learning} the best possible approximation for the reference configuration, has been adopted. The training of NNsPOD lasted considerably longer than the previous linear simulation reported in Section \ref{NNsPOD} with the convergence time towards the threshold amounting to $4$ times what was measured ($27$ hours circa), on the same hardware, for the previous case with the linear setting. 
 
 \begin{figure}[H]
\centering
\begin{tikzpicture}

\definecolor{color2}{rgb}{0,0.75,0}
\definecolor{color1}{rgb}{0.75,0,0.75}
\definecolor{color0}{rgb}{0,0.75,0.75}

\begin{axis}[
axis on top,
legend cell align={left},
legend style={font=\tiny},
tick pos=both,
xtick style={color=black},
ytick style={color=black},
xlabel={Mode},
ymode=log,
grid=both,
legend pos=north east,
grid style={line width=.1pt, draw=gray!10},
ylabel shift = -5pt,
ylabel={Singular value (log-scale)}
]
\addplot table [color=color0, x expr=\coordindex+1, y index=0] {data/Multiphase_SV_POD};
\addlegendentry{POD}
\addplot table [color=color2, x expr=\coordindex+1, y index=0] {data/Multiphase_SV_NNsPOD};
\addlegendentry{Automatic shift-detected POD}
\end{axis}

\end{tikzpicture}
\caption{Comparison of singular values decay for the POD prior and following NNsPOD's pre-processing bijective mapping. We can appreciate a substantial improvement in accuracy retained per number of modes, in line with the shift-based treatment of the snapshot. This result highlights the possibility of NNsPOD to replicate the results of sPOD (with arbitrary accuracy) to those models for which the manual construction of an exact shift operator is not possible.}
\label{fig7}
\end{figure}
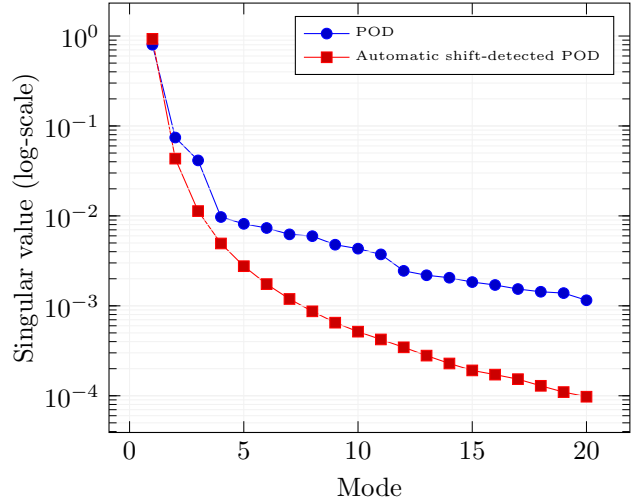
 
 Although this is partially attributed with the higher dimensional full-order discrete space $\mathcal{V}_h$ we did observe a non-trivial behaviour for the convergence of the shift-detection towards the reference configuration. The non-linearity of the FOM in $(\ref{eq10})$ forced the sampling of NNsPOD to not be constrained to simple linear backward shifts which is indeed the goal for the generalisation process itself. In fact, as motivated in Section \ref{NNsPOD} and later shown in Figure \ref{fig6} for the non-uniform linear case, NNsPOD's capability is not limited to linear backward shift for a test snapshot but rather to a non-linear stretching of the shifted coordinates. The search for alternative and non-trivial backward maps, as shown by ShiftNet's outputs in Figure \ref{fig9} and \ref{fig10} entails to our algorithm the property of generating automatically the best possible approximation of the backward shift in the form of a bijective map for the linear manifold reconstruction for unknown transport fields.

\begin{figure}[H]
\includegraphics[width=.48\textwidth]{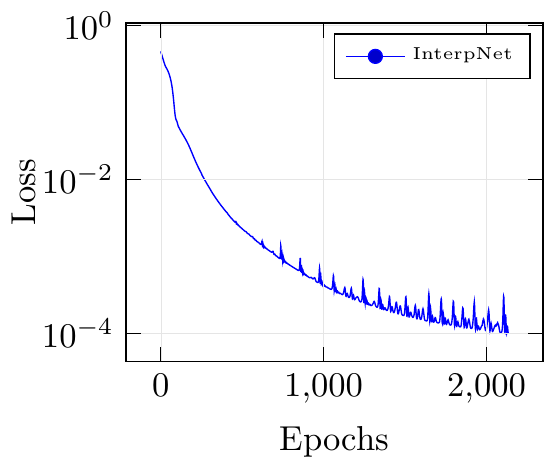}
\includegraphics[width=.48\textwidth]{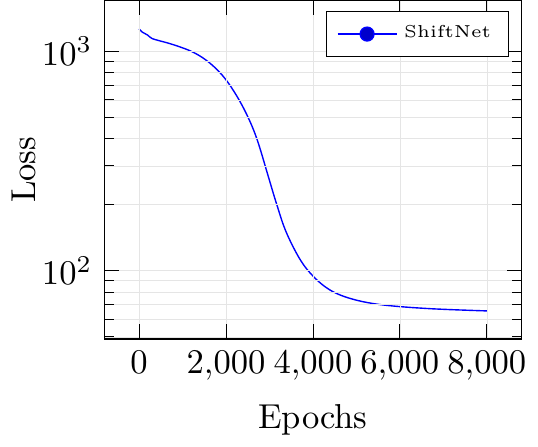}
\caption{NNsPOD training loss optimisation for the multiphase IBVP $(\ref{eq12})$.}
\label{fig12}
\end{figure}

\section{Conclusions}\label{Conclusions}
The purpose of this work was to derive a novel approach for the model order reduction of advection dominated problems. At this aim, the development of NNsPOD allowed us to emphasize the advantages posed by the construction of a non-linear transformation that does not rely on prior knowledge of the physical model. The automatisation of the shift-detection process opposes the sampling of the phase space for the advection field to derive better linear approximations of the manifold at the current state of the art. The multiphase test showcases that the time required for the training of NNsPOD falls within the current estimation for the offline phase with no pre-processing of the snapshots (e.g. traditional POD) while leading to more accurate low-rank subspaces. The online phase of a pre-processing based POD, which is of paramount importance for the industrial applications of any ROM, requires for the backward map to generalise for new instances of the parameter vector; in order to achieve such result the map itself has to be independent on the advection field. To that regard NNsPOD allows for the computation of new solutions during the online phase since its transformation is bijective (given by the continuous dataflow structure of the two networks) and data-driven.

Future works will focus on the validation of the capability of the proposed framework on a more complex setting, e.g. higher-dimensional formulation for the multiphase model in $(\ref{eq10})$ with a parametric treatment of the fluid viscosity.
Moreover, since this work has the main purpose of presenting the NNsPOD methodology for improving the classical POD reduction, we have postponed to future works the integration of such method in a Galerkin projection framework.

Of course, the NNsPOD effectiveness strongly depends on the problem at hand and on the machine learning architecture used for both the networks --- Shiftnet and Interpnet ---.  It is impossible indeed, at the current state, to predict rigorously the optimal structure of the two networks with the choice of the activation function being the sole setting that can be deduced by the FOM solution field as outlined in Sections $\ref{NNsPOD}$ and $\ref{Multiphase}$. Future studies will better investigate the sensitivity of the method with respect to the networks architecture, considering a set of hyper-parameters larger than the one analysed in the current article.

\section*{Acknowledgments}
We acknowledge the support provided by the European Research Council Executive Agency by the Consolidator Grant project AROMA-CFD “Advanced Reduced Order Methods with Applications in Computational Fluid Dynamics” - GA 681447, H2020-ERC CoG 2015 AROMA-CFD, MIUR PRIN NA\_FROM\_PDEs 2019, MIUR FARE-X-AROMA-CFD 2017 and INdAM-GNCS 2019-2020 projects.We acknowledge a long lasting collaboration between Politecnico di Torino (DISMA) and SISSA (mathLab) with exchange of master degree students in Math. Engineering and we are grateful to Prof. Claudio Canuto for insights and continuous support, as well as research collaboration in the framework of MIUR PRIN project NA-FROM-PDEs and H2020 RISE ARIA project.

\bibliographystyle{elsarticle-num} 
\bibliography{citations.bib}

\begin{thebibliography}{10}
\expandafter\ifx\csname url\endcsname\relax
  \def\url#1{\texttt{#1}}\fi
\expandafter\ifx\csname urlprefix\endcsname\relax\def\urlprefix{URL }\fi
\expandafter\ifx\csname href\endcsname\relax
  \def\href#1#2{#2} \def\path#1{#1}\fi

\bibitem{quarteroni_reduced_2014}
A.~Quarteroni, G.~Rozza (Eds.),
  \href{http://link.springer.com/10.1007/978-3-319-02090-7}{Reduced {Order}
  {Methods} for {Modeling} and {Computational} {Reduction}}, Springer
  International Publishing, Cham, 2014.
\newblock \href {https://doi.org/10.1007/978-3-319-02090-7}
  {\path{doi:10.1007/978-3-319-02090-7}}.
\newline\urlprefix\url{http://link.springer.com/10.1007/978-3-319-02090-7}

\bibitem{stein_model_2017}
F.~Chinesta, A.~Huerta, G.~Rozza, K.~Willcox,
  \href{http://doi.wiley.com/10.1002/9781119176817.ecm2110}{Model {Reduction}
  {Methods}}, in: E.~Stein, R.~de~Borst, T.~J.~R. Hughes (Eds.), Encyclopedia
  of {Computational} {Mechanics} {Second} {Edition}, John Wiley \& Sons, Ltd,
  Chichester, UK, 2017, pp. 1--36.
\newblock \href {https://doi.org/10.1002/9781119176817.ecm2110}
  {\path{doi:10.1002/9781119176817.ecm2110}}.
\newline\urlprefix\url{http://doi.wiley.com/10.1002/9781119176817.ecm2110}

\bibitem{quarteroni_reduced_2016}
A.~Quarteroni, A.~Manzoni, F.~Negri,
  \href{https://www.springer.com/gp/book/9783319154305}{Reduced {Basis}
  {Methods} for {Partial} {Differential} {Equations}: {An} {Introduction}}, La
  {Matematica} per il 3+2, Springer International Publishing, 2016.
\newblock \href {https://doi.org/10.1007/978-3-319-15431-2}
  {\path{doi:10.1007/978-3-319-15431-2}}.
\newline\urlprefix\url{https://www.springer.com/gp/book/9783319154305}

\bibitem{hesthaven_certified_2016}
J.~S. Hesthaven, G.~Rozza, B.~Stamm,
  \href{https://www.springer.com/gp/book/9783319224695}{Certified {Reduced}
  {Basis} {Methods} for {Parametrized} {Partial} {Differential} {Equations}},
  {SpringerBriefs} in {Mathematics}, Springer International Publishing, 2016.
\newblock \href {https://doi.org/10.1007/978-3-319-22470-1}
  {\path{doi:10.1007/978-3-319-22470-1}}.
\newline\urlprefix\url{https://www.springer.com/gp/book/9783319224695}

\bibitem{stabile_pod-galerkin_2017}
G.~Stabile, S.~Hijazi, A.~Mola, S.~Lorenzi, G.~Rozza,
  \href{http://arxiv.org/abs/1701.03424}{Pod-{Galerkin} {Reduced} {Order}
  {Methods} for {CFD} {Using} {Finite} {Volume} {Discretisation}: {Vortex}
  {Shedding} {Around} a {Circular} {Cylinder}}, Communications in Applied and
  Industrial Mathematics 8~(1) (2017) 210--236, arXiv: 1701.03424.
\newblock \href {https://doi.org/10.1515/caim-2017-0011}
  {\path{doi:10.1515/caim-2017-0011}}.
\newline\urlprefix\url{http://arxiv.org/abs/1701.03424}

\bibitem{girfoglio_pod-galerkin_2020}
M.~Girfoglio, A.~Quaini, G.~Rozza, \href{http://arxiv.org/abs/2009.13593}{A
  {POD}-{Galerkin} reduced order model for a {LES} filtering approach}, Journal
  of Computational Physics 436 (2021) 110260, arXiv: 2009.13593.
\newblock \href {https://doi.org/10.1016/j.jcp.2021.110260}
  {\path{doi:10.1016/j.jcp.2021.110260}}.
\newline\urlprefix\url{http://arxiv.org/abs/2009.13593}

\bibitem{ballarin_supremizer_2015}
F.~Ballarin, A.~Manzoni, A.~Quarteroni, G.~Rozza,
  \href{http://doi.wiley.com/10.1002/nme.4772}{Supremizer stabilization of
  {POD}-{Galerkin} approximation of parametrized steady incompressible
  {Navier}-{Stokes} equations}, International Journal for Numerical Methods in
  Engineering 102~(5) (2015) 1136--1161.
\newblock \href {https://doi.org/10.1002/nme.4772}
  {\path{doi:10.1002/nme.4772}}.
\newline\urlprefix\url{http://doi.wiley.com/10.1002/nme.4772}

\bibitem{ballarin_pod-galerkin_2016}
F.~Ballarin, G.~Rozza,
  \href{http://doi.wiley.com/10.1002/fld.4252}{{POD}-{Galerkin} monolithic
  reduced order models for parametrized fluid-structure interaction problems:},
  International Journal for Numerical Methods in Fluids 82~(12) (2016)
  1010--1034.
\newblock \href {https://doi.org/10.1002/fld.4252}
  {\path{doi:10.1002/fld.4252}}.
\newline\urlprefix\url{http://doi.wiley.com/10.1002/fld.4252}

\bibitem{hijazi_data-driven_2020}
S.~Hijazi, G.~Stabile, A.~Mola, G.~Rozza,
  \href{http://arxiv.org/abs/1907.09909}{Data-{Driven} {POD}-{Galerkin}
  {Reduced} {Order} {Model} for {Turbulent} {Flows}}, Journal of Computational
  Physics 416 (2020) 109513, arXiv: 1907.09909.
\newblock \href {https://doi.org/10.1016/j.jcp.2020.109513}
  {\path{doi:10.1016/j.jcp.2020.109513}}.
\newline\urlprefix\url{http://arxiv.org/abs/1907.09909}

\bibitem{cohen_kolmogorov_2015}
A.~Cohen, R.~DeVore, \href{https://doi.org/10.1093/imanum/dru066}{{Kolmogorov
  widths under holomorphic mappings}}, IMA Journal of Numerical Analysis 36~(1)
  (2015) 1--12.
\newblock \href {https://doi.org/10.1093/imanum/dru066}
  {\path{doi:10.1093/imanum/dru066}}.
\newline\urlprefix\url{https://doi.org/10.1093/imanum/dru066}

\bibitem{nonino_overcoming_2019}
M.~Nonino, F.~Ballarin, G.~Rozza, Y.~Maday,
  \href{http://arxiv.org/abs/1911.06598}{Overcoming slowly decaying
  {Kolmogorov} n-width by transport maps: application to model order reduction
  of fluid dynamics and fluid--structure interaction problems},
  arXiv:1911.06598 [cs, math]ArXiv: 1911.06598 (Nov. 2019).
\newline\urlprefix\url{http://arxiv.org/abs/1911.06598}

\bibitem{iollo_advection_2014}
A.~Iollo, D.~Lombardi,
  \href{https://link.aps.org/doi/10.1103/PhysRevE.89.022923}{Advection modes by
  optimal mass transfer}, Physical Review E 89~(2) (2014) 022923.
\newblock \href {https://doi.org/10.1103/PhysRevE.89.022923}
  {\path{doi:10.1103/PhysRevE.89.022923}}.
\newline\urlprefix\url{https://link.aps.org/doi/10.1103/PhysRevE.89.022923}

\bibitem{pacciarini_reduced_2015}
P.~Pacciarini, G.~Rozza, Reduced {Basis} {Approximation} of {Parametrized}
  {Advection}-{Diffusion} {PDEs} with {High} {Péclet} {Number}, in:
  A.~Abdulle, S.~Deparis, D.~Kressner, F.~Nobile, M.~Picasso (Eds.), Numerical
  {Mathematics} and {Advanced} {Applications} - {ENUMATH} 2013, Lecture {Notes}
  in {Computational} {Science} and {Engineering}, Springer International
  Publishing, Cham, 2015, pp. 419--426.
\newblock \href {https://doi.org/10.1007/978-3-319-10705-9_41}
  {\path{doi:10.1007/978-3-319-10705-9_41}}.

\bibitem{torlo_stabilized_2018}
D.~Torlo, F.~Ballarin, G.~Rozza,
  \href{http://arxiv.org/abs/1711.11275}{Stabilized weighted reduced basis
  methods for parametrized advection dominated problems with random inputs},
  SIAM/ASA Journal on Uncertainty Quantification 6~(4) (2018) 1475--1502,
  arXiv: 1711.11275.
\newblock \href {https://doi.org/10.1137/17M1163517}
  {\path{doi:10.1137/17M1163517}}.
\newline\urlprefix\url{http://arxiv.org/abs/1711.11275}

\bibitem{rim_model_2018}
D.~Rim, K.~T. Mandli, \href{http://arxiv.org/abs/1805.05938}{Model reduction of
  a parametrized scalar hyperbolic conservation law using displacement
  interpolation}, arXiv:1805.05938 [math]ArXiv: 1805.05938 (May 2018).
\newline\urlprefix\url{http://arxiv.org/abs/1805.05938}

\bibitem{rim_transport_2018}
D.~Rim, S.~Moe, R.~J. LeVeque,
  \href{https://epubs.siam.org/doi/10.1137/17M1113679}{Transport {Reversal} for
  {Model} {Reduction} of {Hyperbolic} {Partial} {Differential} {Equations}},
  SIAM/ASA Journal on Uncertainty Quantification 6~(1) (2018) 118--150.
\newblock \href {https://doi.org/10.1137/17M1113679}
  {\path{doi:10.1137/17M1113679}}.
\newline\urlprefix\url{https://epubs.siam.org/doi/10.1137/17M1113679}

\bibitem{chetverushkin_model_2019}
N.~Cagniart, Y.~Maday, B.~Stamm,
  \href{http://link.springer.com/10.1007/978-3-319-78325-3_10}{Model {Order}
  {Reduction} for {Problems} with {Large} {Convection} {Effects}}, in: B.~N.
  Chetverushkin, W.~Fitzgibbon, Y.~Kuznetsov, P.~Neittaanmäki, J.~Periaux,
  O.~Pironneau (Eds.), Contributions to {Partial} {Differential} {Equations}
  and {Applications}, Vol.~47, Springer International Publishing, Cham, 2019,
  pp. 131--150.
\newblock \href {https://doi.org/10.1007/978-3-319-78325-3_10}
  {\path{doi:10.1007/978-3-319-78325-3_10}}.
\newline\urlprefix\url{http://link.springer.com/10.1007/978-3-319-78325-3_10}

\bibitem{nair_transported_2019}
N.~J. Nair, M.~Balajewicz,
  \href{http://doi.wiley.com/10.1002/nme.5998}{Transported snapshot model order
  reduction approach for parametric, steady-state fluid flows containing
  parameter-dependent shocks: {Model} order reduction for fluid flows
  containing shocks}, International Journal for Numerical Methods in
  Engineering 117~(12) (2019) 1234--1262.
\newblock \href {https://doi.org/10.1002/nme.5998}
  {\path{doi:10.1002/nme.5998}}.
\newline\urlprefix\url{http://doi.wiley.com/10.1002/nme.5998}

\bibitem{taddei_registration_2020}
T.~Taddei, \href{https://epubs.siam.org/doi/10.1137/19M1271270}{A
  {Registration} {Method} for {Model} {Order} {Reduction}: {Data} {Compression}
  and {Geometry} {Reduction}}, SIAM Journal on Scientific Computing 42~(2)
  (2020) A997--A1027.
\newblock \href {https://doi.org/10.1137/19M1271270}
  {\path{doi:10.1137/19M1271270}}.
\newline\urlprefix\url{https://epubs.siam.org/doi/10.1137/19M1271270}

\bibitem{kashima_non-linear_2016}
K.~Kashima, Nonlinear model reduction by deep autoencoder of noise response
  data, 2016 IEEE 55th Conference on Decision and Control (CDC) (2016)
  5750--5755.

\bibitem{hartman_deep-learning_2017}
D.~{Hartman}, L.~K. {Mestha}, A deep learning framework for model reduction of
  dynamical systems, in: 2017 IEEE Conference on Control Technology and
  Applications (CCTA), 2017, pp. 1917--1922.
\newblock \href {https://doi.org/10.1109/CCTA.2017.8062736}
  {\path{doi:10.1109/CCTA.2017.8062736}}.

\bibitem{crisovan_model_2018}
R.~Crisovan, D.~Torlo, R.~Abgrall, S.~Tokareva, Model order reduction for
  parametrized nonlinear hyperbolic problems as an application to uncertainty
  quantification, Journal of Computational and Applied Mathematics 348 (2019)
  466--489.
\newblock \href {https://doi.org/10.1016/j.cam.2018.09.018}
  {\path{doi:10.1016/j.cam.2018.09.018}}.

\bibitem{hoang_projection_2021}
C.~Hoang, K.~Chowdhary, K.~Lee, J.~Ray, Projection-based model reduction of
  dynamical systems using space-time subspace and machine learning (2021).
\newblock \href {http://arxiv.org/abs/2102.03505} {\path{arXiv:2102.03505}}.

\bibitem{reiss_shifted_2018}
J.~Reiss, P.~Schulze, J.~Sesterhenn, V.~Mehrmann,
  \href{https://epubs.siam.org/doi/10.1137/17M1140571}{The {Shifted} {Proper}
  {Orthogonal} {Decomposition}: {A} {Mode} {Decomposition} for {Multiple}
  {Transport} {Phenomena}}, SIAM Journal on Scientific Computing 40~(3) (2018)
  A1322--A1344.
\newblock \href {https://doi.org/10.1137/17M1140571}
  {\path{doi:10.1137/17M1140571}}.
\newline\urlprefix\url{https://epubs.siam.org/doi/10.1137/17M1140571}

\bibitem{sarna_hyper-reduction_2021}
N.~Sarna, S.~Grundel, \href{http://arxiv.org/abs/2003.06362}{Hyper-reduction
  for parametrized transport dominated problems via online-adaptive reduced
  meshes}, arXiv:2003.06362 [cs, math]ArXiv: 2003.06362 (Jan. 2021).
\newline\urlprefix\url{http://arxiv.org/abs/2003.06362}

\bibitem{lee_model_2019}
K.~Lee, K.~T. Carlberg,
  \href{https://www.sciencedirect.com/science/article/pii/S0021999119306783}{Model
  reduction of dynamical systems on nonlinear manifolds using deep
  convolutional autoencoders}, Journal of Computational Physics 404 (2020)
  108973.
\newblock \href {https://doi.org/https://doi.org/10.1016/j.jcp.2019.108973}
  {\path{doi:https://doi.org/10.1016/j.jcp.2019.108973}}.
\newline\urlprefix\url{https://www.sciencedirect.com/science/article/pii/S0021999119306783}

\bibitem{torlo_model_2020}
D.~Torlo, \href{http://arxiv.org/abs/2003.13735}{Model {Reduction} for
  {Advection} {Dominated} {Hyperbolic} {Problems} in an {ALE} {Framework}:
  {Offline} and {Online} {Phases}}, arXiv:2003.13735 [cs, math]ArXiv:
  2003.13735 (Mar. 2020).
\newline\urlprefix\url{http://arxiv.org/abs/2003.13735}

\bibitem{mojgani_arbitrary_2017}
R.~Mojgani, M.~Balajewicz,
  \href{http://adsabs.harvard.edu/abs/2017APS..DFD.M1008M}{Arbitrary
  {Lagrangian} {Eulerian} framework for efficient projection-based reduction of
  convection dominated nonlinear flows} (2017) M1.008.
\newline\urlprefix\url{http://adsabs.harvard.edu/abs/2017APS..DFD.M1008M}

\bibitem{peng_learning-based_2021}
Z.~Peng, M.~Wang, F.~Li, \href{http://arxiv.org/abs/2105.14633}{A
  learning-based projection method for model order reduction of transport
  problems}, arXiv:2105.14633 [cs, math]ArXiv: 2105.14633 (May 2021).
\newline\urlprefix\url{http://arxiv.org/abs/2105.14633}

\bibitem{moukalled_finite_2016}
F.~Moukalled, L.~Mangani, M.~Darwish, The {Finite} {Volume} {Method} in
  {Computational} {Fluid} {Dynamics}: {An} {Advanced} {Introduction} with
  {OpenFOAM}® and {Matlab}, 1st Edition, no. 113 in Fluid {Mechanics} and
  {Its} {Applications}, Springer International Publishing : Imprint: Springer,
  Cham, 2016.

\bibitem{atkinson_num_anal}
K.~Atkinson, W.~Han, Theoretical Numerical Analysis: A Functional Analysis
  Framework, Vol.~39, 2009.
\newblock \href {https://doi.org/10.1007/978-1-4419-0458-4}
  {\path{doi:10.1007/978-1-4419-0458-4}}.

\bibitem{quarteroni_numerical_2009}
A.~Quarteroni, Numerical models for differential problems, no. v. 2 in {MS} \&
  {A}, Springer, Milan ; New York, 2009, oCLC: ocn288986457.

\bibitem{leonard_quick_1979}
B.~Leonard,
  \href{https://www.sciencedirect.com/science/article/pii/0045782579900343}{A
  stable and accurate convective modelling procedure based on quadratic
  upstream interpolation}, Computer Methods in Applied Mechanics and
  Engineering 19~(1) (1979) 59--98.
\newblock \href {https://doi.org/https://doi.org/10.1016/0045-7825(79)90034-3}
  {\path{doi:https://doi.org/10.1016/0045-7825(79)90034-3}}.
\newline\urlprefix\url{https://www.sciencedirect.com/science/article/pii/0045782579900343}

\bibitem{stewart_svd_1993}
G.~W. Stewart, \href{https://doi.org/10.1137/1035134}{On the early history of
  the singular value decomposition}, SIAM Review 35~(4) (1993) 551--566.
\newblock \href {http://arxiv.org/abs/https://doi.org/10.1137/1035134}
  {\path{arXiv:https://doi.org/10.1137/1035134}}, \href
  {https://doi.org/10.1137/1035134} {\path{doi:10.1137/1035134}}.
\newline\urlprefix\url{https://doi.org/10.1137/1035134}

\bibitem{peherstorfer_model_2020}
B.~Peherstorfer, \href{https://doi.org/10.1137/19M1257275}{Model reduction for
  transport-dominated problems via online adaptive bases and adaptive
  sampling}, SIAM Journal on Scientific Computing 42~(5) (2020) A2803--A2836.
\newblock \href {http://arxiv.org/abs/https://doi.org/10.1137/19M1257275}
  {\path{arXiv:https://doi.org/10.1137/19M1257275}}, \href
  {https://doi.org/10.1137/19M1257275} {\path{doi:10.1137/19M1257275}}.
\newline\urlprefix\url{https://doi.org/10.1137/19M1257275}

\bibitem{frank_machine-learning_2020}
M.~Frank, D.~Drikakis, V.~Charissis,
  \href{https://www.mdpi.com/2079-3197/8/1/15}{Machine-{Learning} {Methods} for
  {Computational} {Science} and {Engineering}}, Computation 8~(1) (2020) 15.
\newblock \href {https://doi.org/10.3390/computation8010015}
  {\path{doi:10.3390/computation8010015}}.
\newline\urlprefix\url{https://www.mdpi.com/2079-3197/8/1/15}

\bibitem{kutz_deep_2017}
J.~N. Kutz,
  \href{https://www.cambridge.org/core/product/identifier/S002211201600803X/type/journal_article}{Deep
  learning in fluid dynamics}, Journal of Fluid Mechanics 814 (2017) 1--4.
\newblock \href {https://doi.org/10.1017/jfm.2016.803}
  {\path{doi:10.1017/jfm.2016.803}}.
\newline\urlprefix\url{https://www.cambridge.org/core/product/identifier/S002211201600803X/type/journal_article}

\bibitem{sarghini_neural_2003}
F.~Sarghini, G.~de~Felice, S.~Santini,
  \href{https://linkinghub.elsevier.com/retrieve/pii/S0045793001000986}{Neural
  networks based subgrid scale modeling in large eddy simulations}, Computers
  \& Fluids 32~(1) (2003) 97--108.
\newblock \href {https://doi.org/10.1016/S0045-7930(01)00098-6}
  {\path{doi:10.1016/S0045-7930(01)00098-6}}.
\newline\urlprefix\url{https://linkinghub.elsevier.com/retrieve/pii/S0045793001000986}

\bibitem{renardy_numerical_2001}
M.~Renardy, Y.~Renardy, J.~Li,
  \href{https://linkinghub.elsevier.com/retrieve/pii/S0021999101967853}{Numerical
  {Simulation} of {Moving} {Contact} {Line} {Problems} {Using} a
  {Volume}-of-{Fluid} {Method}}, Journal of Computational Physics 171~(1)
  (2001) 243--263.
\newblock \href {https://doi.org/10.1006/jcph.2001.6785}
  {\path{doi:10.1006/jcph.2001.6785}}.
\newline\urlprefix\url{https://linkinghub.elsevier.com/retrieve/pii/S0021999101967853}

\bibitem{okagaki_numerical_2021}
Y.~Okagaki, T.~Yonomoto, M.~Ishigaki, Y.~Hirose,
  \href{https://www.mdpi.com/2311-5521/6/2/80}{Numerical {Study} on an
  {Interface} {Compression} {Method} for the {Volume} of {Fluid} {Approach}},
  Fluids 6~(2) (2021) 80.
\newblock \href {https://doi.org/10.3390/fluids6020080}
  {\path{doi:10.3390/fluids6020080}}.
\newline\urlprefix\url{https://www.mdpi.com/2311-5521/6/2/80}

\end{thebibliography}

\end{document}